\documentclass[11pt]{amsart}
\usepackage{amsmath, amssymb, amscd, mathrsfs, slashed, enumerate, url}
\usepackage{amsfonts}
\usepackage{color}
\usepackage{extsizes}
\usepackage{xcolor}
\usepackage[all,cmtip]{xy}
\usepackage{multicol}
\usepackage{enumitem} 
\usepackage{indentfirst}
\usepackage{latexsym}
\usepackage{bm}
\usepackage{graphicx}
\usepackage{subfigure,esint}
\usepackage{float}
\usepackage{verbatim}
\usepackage{comment}
\usepackage{esint}
\usepackage[top=1in, bottom=1in, left=1.25in, right=1.25in]{geometry}
\usepackage{epsfig,dsfont,amsthm,amsfonts,amsbsy}
\usepackage[colorlinks]{hyperref}
\usepackage[toc,page]{appendix}
\usepackage{tikz-cd}
\usepackage{amsthm,xcolor}
\newtheorem{theorem}{Theorem}[section]

\newtheorem{conjecture}[theorem]{Conjecture}
\newtheorem{corollary}[theorem]{Corollary}
\newtheorem{proposition}[theorem]{Proposition}
\newtheorem{definition}[theorem]{Definition}
\newtheorem{lemma}[theorem]{Lemma}

\theoremstyle{remark}

\newtheorem{example}[theorem]{Example}

\newtheorem{remark}[theorem]{Remark}

\newcommand{\lan}{\langle }
\newcommand{\ran}{\rangle}

\newcommand{\tM}{\widetilde{M}}
\newcommand{\SU}{\mathrm{SU}}

\newcommand{\MH}{\mathcal{H}}

\newcommand{\MM}{\mathcal{M}}

\newcommand{\MC}{\mathcal{C}}
\newcommand{\Tr}{\mathrm{Tr}}

\newcommand{\Ric}{\mathrm{Ric}}

\newcommand{\bu}{\mathbf{u}}
\newcommand{\bv}{\mathbf{v}}

\newcommand{\MT}{\mathcal{T}}

\newcommand{\ZT}{\mathbb{Z}/2}

\newcommand{\vol}{\mathrm{vol}}

\newcommand{\id}{\mathrm{id}}

\newcommand{\BXi}{\overline{\Xi}}

\newcommand{\SLC}{\mathrm{SL}_2(\mathbb{C})}

\newcommand{\SL}{\mathrm{SL}}

\newcommand{\inj}{\mathrm{inj}}

\newcommand{\SO}{\mathrm{SO}}

\newcommand{\mbR}{\mathbb{R}}

\newcommand{\mbM}{\mathbb{M}}

\let\Re\undefined
\newcommand{\Re}{\mathrm{Re}}

\newcommand{\cX}{\mathcal{X}}

\newcommand{\BcX}{\overline{\cX(\Gamma)}}

\newcommand{\MPL}{\mathcal{PL}}
\newcommand{\MMZ}{\mathcal{M}_{\ZT}}
\newcommand{\MMc}{\mathcal{M}_{\mathrm{SL}_2(\mathbb{C})}}
\newcommand{\BMMc}{\overline{\mathcal{M}}_{\mathrm{SL}_2(\mathbb{C})}}
\newcommand{\tbv}{\tilde{\bv}}

\usepackage{calligra}
\DeclareMathAlphabet{\mathcalligra}{T1}{calligra}{m}{n}

\theoremstyle{plain}
\makeatletter
\newtheorem*{rep@theorem}{\rep@title}
\newcommand{\newreptheorem}[2]{%
	\newenvironment{rep#1}[1]{%
		\def\rep@title{#2 \ref{##1}}%
		\begin{rep@theorem}}%
		{\end{rep@theorem}}}
\makeatother
\newreptheorem{theorem}{Theorem}
\newreptheorem{proposition}{Theorem}

\usepackage[utf8]{inputenc}
\usepackage[english]{babel}
\usepackage{fancyhdr}
\usepackage{accents}

\begin{document}
	\title[Harmonic 1-forms and the Morgan-Shalen compactification]{$\ZT$ harmonic 1-forms,
 $\mbR$-trees, and\\ the   Morgan-Shalen compactification}
\author[He]{Siqi He}
\email{sqhe@amss.ac.cn}
\address{Morningside Center of Mathematics,
	Chinese Academy of Sciences, 
	Beijing, 100190 China}

\author[Wentworth]{Richard  Wentworth}
\email{raw@umd.edu}
\address{Department of Mathematics,
	University of Maryland,
	College Park, MD 20742, USA
}

\author[Zhang]{Boyu Zhang}
\email{bzh@umd.edu}
\address{
	Department of Mathematics,
	University of Maryland,
	College Park, MD 20742, USA
}

\makeatletter
\@namedef{subjclassname@2020}{%
  \textup{2020} Mathematics Subject Classification}
\makeatother

\date{\today}

\begin{abstract}
This paper studies the relationship between an analytic compactification of the
   moduli space of flat  $\mathrm{SL}_2(\mathbb{C})$
connections on a closed, oriented 
 3-manifold $M$  defined by Taubes, and the Morgan--Shalen compactification
 of the $\mathrm{SL}_2(\mathbb{C})$ character variety of the fundamental
    group of $M$. We exhibit an explicit  correspondence between $\mathbb{Z}/2$ harmonic 1-forms,
    measured foliations, and equivariant harmonic maps to $\mathbb{R}$-trees,
    as initially proposed by Taubes. As an application, we prove that $\mathbb{Z}/2$ harmonic 1-forms exist on all Haken manifolds with respect to all Riemannian metrics. We also show that there exist manifolds   that support singular  $\mathbb{Z}/2$ harmonic 1-forms but have compact $\SL_2(\mathbb{C})$ character varieties, which resolves a folklore conjecture. 
\end{abstract}
	\maketitle

\section{Introduction}
Let $(M, g)$ be a closed, oriented Riemannian 3-manifold with fundamental group
$\Gamma = \pi_1(M)$. Let $\cX(\Gamma)$ denote the $\SLC$ character variety
of $\Gamma$ as defined in, for example, \cite{cullershalen1983varieties}. 
Note that as a set, the closed points of $\cX(\Gamma)$ correspond to conjugacy  classes of
\emph{completely reducible representations}, i.e.\ those that are either
irreducible or direct sums of a $1$--dimensional representation and its
dual. The space $\cX(\Gamma)$ has the structure of an affine
algebraic variety which, if positive dimensional, admits an ${\rm Aut}(\Gamma)$--invariant
compactification $\BcX$, called the Morgan--Shalen compactification 
\cite{cullershalen1983varieties, morganshalen1984valuationstrees, morganshalen88degeneration, morganshalen88degernerationtwo}.
A boundary point in $\BcX$ is given by the projective class of a  
length function for an isometric action of $\Gamma$ on an $\mathbb{R}$--tree.

In \cite{Taubes20133manifoldcompactness}, Taubes introduced a compactification of the space of
flat $\SLC$ connections on $3$--manifolds from an analytic perspective.
Consider the set of solutions to the equations
\begin{equation}
	\label{eqn_KW_system}
	F_A = \frac12 [\phi,\phi],\quad d_A\phi = 0, \quad  d_A^*\phi = 0,
\end{equation}
where $A$ is a connection on an $\SU(2)$ bundle $P$ over $M$, and $\phi$ is
a section of $T^*M\otimes \mathfrak{g}_P$, where $\mathfrak{g}_P=
P\times_{ad}\mathfrak{su}(2)$ is the adjoint bundle. 
Let $\MMc$ denote the moduli space of solutions to \eqref{eqn_KW_system}
 up to $\SU(2)$ gauge transformations. 
We abuse notation and also use $\cX(\Gamma)$ to denote the topological space of closed points of $\cX(\Gamma)$ with the analytic topology.
Then by the results in \cite{donaldson1987twisted} and
\cite[Thm.\ 3.3]{corlette1988flat}, $\MMc$ is canonically homeomorphic to
$\cX(\Gamma)$.
The compactness results of Taubes \cite{Taubes20133manifoldcompactness} define a compactification of $\MMc$, where the boundary points are described by a class of objects called \emph{$\mathbb{Z}/2$ harmonic $1$--forms}.
A $\mathbb{Z}/2$ harmonic $1$--form can be regarded as a generalization to
$3$-manifolds of holomorphic  quadratic differentials on a Riemann surface. 
Its definition and properties will be reviewed in Section \ref{sec_Z2_forms} below.

This article studies the relationship between the Morgan--Shalen compactification and the Taubes compactification. 
The Morgan--Shalen compactification is closely related to topological concepts such as singular measured foliations,
incompressible surfaces, and $\mathbb{R}$--trees.
Our main results will provide a topological interpretation of the analytic limits defined by Taubes. 
The relationship between the Morgan--Shalen compactification and the analytical compactification for 
two-dimensional manifolds was studied in \cite{daskalopoulos2000morganshalen}; 
see also \cite{katzarkov2015harmonic,he2023algebraic, ott2020higgs,loftin2022limit,burger2023real} for extensions. 

Our results are based on the theory of harmonic maps from manifolds to
$\mathbb{R}$--trees and more generally to metric spaces of nonpositive
curvature (\emph{NPC spaces}) as developed by Wolf and   Korevaar--Schoen
\cite{wolf1995harmonic, korevaar1993sobolev, korevaari1997global}. 
Let $\widetilde{M}$ be the universal cover of $M$. 
By the Corlette--Donaldson theorem, every completely reducible flat $\SLC$
connection on $M$ defines a $\pi_1(M)$--equivariant
harmonic map from $\widetilde{M}$ to the $3$--dimensional hyperbolic space form $\mathbb{H}^3$. 
The Gromov--Hausdorff limit of the convex hulls of the images of a divergent
sequence of harmonic maps is the 
corresponding Morgan--Shalen tree (cf.\ \cite{daskalopoulos1998character}).
We establish a direct relationship between the $\ZT$ harmonic $1$--forms
appearing in the  Taubes limit and the limiting harmonic maps to $\mathbb{R}$--trees.
Such a relationship was proposed by Taubes in \cite[pp.\ 12-14]{taubes2013compactness}.

More precisely, let $\MMZ$ be the moduli space of $\ZT$ harmonic 1-forms up
to rescaling, equipped  with  the $L^2$ topology (see \eqref{eqn_MMZ_def}),
and
let $\BMMc$ denote the Taubes compactification with boundary $\partial \BMMc \subset \MMZ$ 
(see Definition \ref{def_BMMc}). Let $\MPL(\Gamma)$ be the space of
projective length functions of minimal $\Gamma$--actions on complete 
$\mathbb{R}$--trees (see Section \ref{subsec_MS_cpt}). 
In Section \ref{subsection_construction_boundary_map}, we construct a map
\begin{equation}
	\label{eq_boundarymap}
	\MH: \MPL(\Gamma) \to \MMZ.
\end{equation}
Using this, we define a map
\begin{equation}
	\label{eq_compactified_map}
	\begin{split}
		\BXi: \BcX\to \BMMc,
	\end{split}
\end{equation}
where $\BXi$ is given by the Riemann--Hilbert correspondence on $\cX(\Gamma)$ and  by $\MH$ on $\partial \BcX$. 
 Then we have the following result.

\begin{theorem}
	\label{thm_map_on_closure}
	The map $\BXi: \BcX \to \BMMc$ is continuous and surjective.
\end{theorem}

Let $\bv \in \partial \BMMc$ be a $\ZT$ harmonic form in the boundary of the compactified moduli space.
In Section \ref{sec_leaf_space}, we show that the leaf space of the pull-back of 
$\bv$ to $\widetilde{M}$ naturally has the structure of an $\mathbb{R}$--tree, 
which we denote by $\MT_{\widetilde{M},\tilde\bv}$. Let $\pi_{\tilde\bv}$ denote the projection from $\widetilde{M}$ to  $\MT_{\widetilde{M},\tilde\bv}$. The following result describes the relationship between the leaf space and the Morgan--Shalen tree.

\begin{theorem}
	\label{thm_factoriztion_tree_maps}
	Let $\bv \in \partial \BMMc$ be a $\ZT$ harmonic form in the boundary of the compactified moduli space, and assume $\ell \in \MPL(\Gamma)$ satisfies $\MH(\ell) = \bv$. Suppose $\MT_\ell$ is the minimal tree with length function $\ell$ and $u:\widetilde{M}\to \MT_\ell$ is the equivariant harmonic map given by the Morgan--Shalen compactification. Assume $u$ is scaled so that its energy is equal to $1$. Then 
	\begin{enumerate}
		\item 	The map $u$ factors uniquely as the composition of $\pi_{\tilde \bv}$ and a continuous map $f: \MT_{\widetilde{M},\tilde\bv}  \to  \MT_\ell$.
		\item The map $f$ is $1$--Lipschitz.
		\item The map $\pi_{\tilde\bv}$ is harmonic. 
		\item $|\nabla u| = |\nabla \pi_{\tilde\bv}|$ (see Section \ref{subsec_energy_density} for the definitions). Away from the zero set of $|\nabla u|$, we have $\nabla u = \nabla\pi_{\tilde\bv}$ as two-valued differential forms.
	\end{enumerate}
\end{theorem}

The above results establish a bridge between the analytic and algebraic compactifications of the moduli space of 
flat $\SLC$ connections on $3$--manifolds. This allows us to prove new
results on one side using results from the other side. We discuss several
such applications in Section \ref{sec_applications}.

In Section \ref{subsec_existenceHakenmanifold}, we prove an existence result of $\mathbb{Z}/2$ harmonic forms by considering harmonic maps to $\mathbb{R}$--trees. 

A non-zero $\mathbb{Z}/2$ harmonic 1-form $\bv$ is called \emph{non-trivial} if
$\bv$ has non-trivial holonomy; in other words, if it is not given by a single-valued $1$--form.
We say that $\bv$  is
\emph{singular} if there exists a point on $M$ where $\bv$ cannot be
locally lifted to a  single-valued form.  Clearly, singular implies
non-trivial, but the converse is false in general. 
Using the map $\mathcal{H}$ constructed in Section \ref{subsection_construction_boundary_map}, we prove the following theorem.

\begin{theorem}
	\label{thm_pi1_inj_Z2_form}
	Suppose $M$ is a rational homology sphere, and suppose there exists a
    closed connected embedded surface  $S\subset M$ such that $S$ is two-sided,
    $\pi_1$--injective, and does not bound an embedded ball. Then for every
    Riemannian metric $g$ on $M$, there exists a non-trivial $\mathbb{Z}/2$ harmonic 1-form on $(M,g)$.
\end{theorem}

We remark that while all previous analytical constructions of non-trivial 
$\mathbb{Z}/2$ harmonic 1-forms and spinors require the metric to take specific forms (\cite{doan2017existence, taubes2020examples, taubes24topological,he2022existence, heparker2024notes, chen2024existence}), Theorem \ref{thm_pi1_inj_Z2_form} holds for \emph{all} metrics on the given manifolds. 

A rational homology sphere $M$ satisfies the conditions in Theorem \ref{thm_pi1_inj_Z2_form} if and only if $M$ is reducible or Haken.  
On the other hand, if a closed oriented 3-manifold $M$ has $b_1(M)>0$, then there always exist trivial $\mathbb{Z}/2$ harmonic forms 
that are given by usual
harmonic $1$--forms.  Therefore, Theorem \ref{thm_pi1_inj_Z2_form} has the following consequence.

\begin{corollary}
	\label{cor_haken_exist}
	Suppose a closed oriented 3-manifold $M$ is reducible or Haken. Then for every Riemannian metric $g$ on $M$, there exist $\mathbb{Z}/2$ harmonic 1-forms on $(M, g)$.
\end{corollary}

Since $\mathbb{Z}/2$ harmonic 1-forms describe the boundary of compactified moduli spaces of solutions to \eqref{eqn_KW_system}, it is natural to ask whether all of them can be deformed to solutions of the equation. The Kuranishi structure near the boundary of compactified moduli spaces for generalized Seiberg--Witten equations has been studied by Doan--Walpuski \cite{doan2020deformation} and Parker \cite{parker2024gluing}. Equation \eqref{eqn_KW_system} is a special case of the generalized Seiberg--Witten equations; in general, the compactified moduli spaces are described by $\mathbb{Z}/2$ harmonic \emph{spinors}, which is a generalization of the concept of $\mathbb{Z}/2$ harmonic forms. 

For the Seiberg--Witten equations with two spinors, Parker \cite{parker2024gluing} proved that all $\mathbb{Z}/2$ harmonic spinors, under some nondegenerate conditions, can be realized as limits of 1-parameter families of solutions to the equation.
Motivated by Parker's result, there has been a folklore conjecture that $\mathbb{Z}/2$ harmonic 1-forms could not exists on any 3-manifold with compact $\mathrm{SL}_2(\mathbb{C})$ character variety. A precise formulation of this conjecture was recently stated in \cite[Conjecture 1.14]{heparker2024notes}.


Theorem \ref{thm_pi1_inj_Z2_form} implies that there are counterexamples to this folklore conjecture. In fact, we have the following result.

\begin{corollary}
	\label{cor_haken_dim_0}
	There exist infinitely many closed 3-manifolds $M$ such that the
	$\mathrm{SL}_2(\mathbb{C})$ character variety of $\pi_1(M)$ is compact,
	and $M$ supports a singular $\mathbb{Z}/2$ harmonic 1-form with respect to every Riemannian metric. 
\end{corollary}

\begin{proof}
Boyer--Zhang \cite[Theorem 1.8]{boyer98cullershalen} and Motegi \cite{motegi88haken} constructed infinitely many closed oriented Haken 3-manifolds whose $\mathrm{SL}_2(\mathbb{C})$ character varieties are zero-dimensional. Since character varieties are affine varieties, zero-dimensional character varieties are compact. 
There are infinitely many examples $M$ in \cite{motegi88haken} such that $H_1(M)$ are cyclic groups with odd orders. If $M$ is a rational homology sphere such that $H_1(M)$ has no 2-torsion, every $\mathbb{Z}/2$ harmonic form on $M$ is singular. Hence the result is proved by Theorem \ref{thm_pi1_inj_Z2_form}.
\end{proof}

Corollary \ref{cor_haken_dim_0} gives the first example of $\mathbb{Z}/2$
harmonic 1-forms (and more generally, $\mathbb{Z}/2$ harmonic spinors) that cannot be deformed into solutions of the corresponding gauge-theoretic equations. 
It also suggests that some properties of solutions to the Seiberg--Witten equations with two spinors, as studied in \cite{doan2020deformation,doan2017existence,parker2024gluing}, may not generalize to flat $\SL_2(\mathbb{C})$ connections.

In Sections \ref{subsec_proj_leaf_space} to \ref{subsec_Hub_Mas_map}, we give several more applications of Theorems \ref{thm_map_on_closure} and \ref{thm_factoriztion_tree_maps} and discuss some related results.
Using the regularity results of Parker \cite{parker2023concentrating}, we prove a stronger regularity result for the Morgan--Shalen convergence (Corollary \ref{cor_C^infty_KS_lim}).
Using Theorem \ref{thm_factoriztion_tree_maps}, we show that the $\mathbb{Z}/2$ harmonic forms that arise as limits in Taubes' construction must satisfy certain additional properties (Corollary \ref{cor_condition_MMZ_taubes_lim}). 
The Hubbard--Masur map \cite{hubbard1979quadratic} takes the space of measured foliations on a Riemann surface, modulo certain equivalence relations, to the space of quadratic differentials. Wolf \cite{wolf1995harmonic} showed that this map is a homeomorphism.
Our construction of $\MH$ can be interpreted as a generalization of the Hubbard--Masur map in dimension three, and it associates a canonical measured foliation with every boundary point in the Morgan--Shalen limit (see Section \ref{subsec_Hub_Mas_map}).

Finally, we address a folklore non-existence conjecture for $\mathbb{Z}/2$ harmonic $1$--forms on  manifolds that are diffeomorphic to $S^3$, which is motivated by the relation of $\ZT$ harmonic $1$--forms and the $\SLC$ representation variety. We refer to \cite{heparker2024notes} for more 
discussion about this conjecture.
\begin{conjecture}
	\label{conj:Z2existenceonS3}
	Let $g$ be a Riemannian metric on $S^3$.
 There exists no $\ZT$ harmonic $1$--form on $(S^3,g)$.
\end{conjecture}

Note that if a manifold has positive Ricci curvature, then the
non-existence of $\ZT$ harmonic $1$--forms follows directly from the
Weitzenb\"ock formula (see Eqns.\ \eqref{eqn_Weitzenbock_harmonic_form}, \eqref{eqn_int_laplacian_|v|^2}).  

We approach this conjecture using the relationship between $\ZT$ harmonic
$1$--forms and measured foliations, as developed in Sections
\ref{sec_leaf_space} and \ref{sec_applications} below. Given a $\ZT$ harmonic 1-form $\bv$ on $S^3$, the leaf space of the measured foliation defined by $\bv$ is an $\mathbb{R}$--tree (see Theorem \ref{thm_pseudo-metric-R-tree}). 
Under certain assumptions, we show that the projection map from $S^3$ to
the leaf space satisfies a maximum principal, which leads to a
contradiction. As a consequence, we prove the following partial resolution
of Conjecture \ref{conj:Z2existenceonS3}. 
\begin{theorem}
	\label{thm_ZT_s3}
	There is no $\ZT$ harmonic $1$--form $\bv$ on $(S^3,g)$ satisfying both of
    the following conditions:
	\begin{enumerate}
		\item Every zero point of $\bv$ is cylindrical (see Definition \ref{def_cylindrical_zero}).
		\item For every arc $\gamma$ in $M \setminus Z$ transverse to $\ker \bv$, where $Z$ denotes the zero locus of $\bv$, we have 
		$
		\mu_\bv(\gamma) = d_{M,\bv}(x,y).
		$
	\end{enumerate}
\end{theorem}

Condition (1) is a local regularity requirement for the zero locus of
$\bv$; it is weaker than the usual regularity assumption in the study of
$\ZT$ harmonic $1$--forms as stated
\cite{donaldson21deformation,he23branched, parker2023deformations}. The
notations $\mu_\bv$ and $d_{M,\bv}$ in Condition (2) refer to the
transverse invariant measure and the distance function on the leaf space
defined by $\bv$ (see \eqref{eqn_transverse_measure_from_bv} and
\eqref{eq_lengthfunction_defined_by_Z2harmonicform}). We will show that
Condition (2) always holds for $\ZT$ harmonic $1$--forms that appear on the
boundary of Taubes' compactification (Corollary
\ref{cor_condition_MMZ_taubes_lim}, part (2)). We  note that a completely different approach to Conjecture \ref{conj:Z2existenceonS3} is given in a forthcoming work of Parker \cite{parker24gaugetheoretic}, whch uses the gluing constructions in gauge theory.

\textbf{Acknowledgements.} The authors wish to express their gratitude to many people for their interest and helpful comments. Among them are Nathan Dunfield, Xinghua Gao, 
Andriy Haydys, Zhenkun Li, Yi Liu, Ciprian Manolescu, Rafe Mazzeo, Tomasz Mrowka, Yi Ni, Jean-Pierre Otal, Greg Parker,  Ao Sun, Clifford Taubes, Thomas Walpuski and Mike Wolf.

S.H. is partially supported by NSFC grant No.12288201 and No.2023YFA1010500.
R.W.’s research is supported by NSF grant DMS-2204346. B.Z. is partially supported by NSF grant DMS-2405271 and a travel grant from the Simons Foundation.

\section{The Morgan-Shalen compactification, length functions, and $\mbR$-trees}

In this section, we briefly review relevant results about the Morgan--Shalen compactification and $\mathbb{R}$--trees from \cite{cullershalen1983varieties,morganshalen1984valuationstrees, morganshalen88degeneration, morganshalen88degernerationtwo}. For a more comprehensive introduction, we refer the reader to the survey by Otal \cite{otal15surveycompactification}.

\subsection{$\mathbb{R}$--trees and length functions}
An $\mathbb{R}$--tree is a metric space $(\mathcal{T}, d_{\mathcal{T}})$ such that
every pair of points is connected by a unique arc, 
and every arc is isometric to a closed interval in $\mathbb{R}$ as a subspace of $\mathcal{T}$.

Let $\Gamma$ be a finitely generated group. A \emph{$\Gamma$--tree} is an $\mathbb{R}$--tree $\mathcal{T}$ with an isometric action $\rho: \Gamma \to \text{Isom}(\mathcal{T})$. A $\Gamma$--tree is called \emph{minimal} if there is no proper $\Gamma$--invariant subtree. The \emph{length function} is defined by
\begin{equation}
	\label{eq_lengthfunction}
	\ell_{\rho}(\gamma) := \inf_{x \in \mathcal{T}} d_{\mathcal{T}}(x, \rho(\gamma)x).
\end{equation}
By \cite[Sec.\ (1.3)]{cullermorgan97group},
every $\gamma \in \Gamma$ acts semisimply, i.e., the infimum is realized at some point in $\mathcal{T}$. By \cite[Prop. II.2.15]{morganshalen1984valuationstrees}, $\ell_{\rho}$ is identically zero if and only if $\Gamma$ has a fixed point.

In general, $\ell_\rho(\gamma)$ depends only on the conjugacy class of $\rho$ and the conjugacy class of $\gamma$. So let $C$ be the set of conjugacy classes of $\Gamma$, and let $\mathbb{P}(C) := \mathbb{P}(\mathbb{R}^{C})$ be the real projective space of nonzero functions on $C$. If $\ell_{\rho} \not\equiv 0$, then the class of $\ell_{\rho}$ in $\mathbb{P}(C)$ is called the \emph{projective length function}.

A length function is called \emph{abelian} if it is given by $|\mu(\gamma)|$ for some homomorphism $\mu: \Gamma \to \mathbb{R}$. A \emph{ray} $R$ in an $\mathbb{R}$--tree is the image of an isometric embedding of $[0, +\infty)$. A $\Gamma$--tree is said to have a \emph{fixed end} if there exists a ray $R$ such that for every $\gamma \in \Gamma$, $\gamma(R) \cap R$ is also a ray. We will need the following result.

\begin{theorem}[{\cite[Thm.\ 3.7]{cullermorgan97group}}]
	\label{thm_nonabelian_determinant_tree}
	Assume $\mathcal{T}$ is a minimal $\Gamma$--tree with a non-trivial length function $\ell$. Then $\ell$ is non-abelian if and only if $\Gamma$ acts without fixed ends. Moreover, if $\ell$ is non-abelian and $\mathcal{T}'$ is another minimal $\Gamma$--tree with the same length function, then there exists a unique $\Gamma$-equivariant isometry between $\mathcal{T}$ and $\mathcal{T}'$.
\end{theorem}

When $\ell$ is abelian, the $\Gamma$--tree is not always uniquely
determined by the length function. A counterexample is given in \cite[Ex.\ 3.9]{cullermorgan97group}. 
To simplify terminology, we will usually refer to a ``$\Gamma$--tree'' as a ``tree'' when the $\Gamma$--action is clear from the context.

\subsection{The Morgan--Shalen compactification}
\label{subsec_MS_cpt}
The \emph{character} $\chi_{\rho}: \Gamma \to \mathbb{C}$ of a representation $\rho: \Gamma \to \mathrm{SL}_2(\mathbb{C})$ is defined by $\chi_{\rho}(\gamma) := \mathrm{Tr}(\rho(\gamma))$. Clearly, $\chi_\rho$ depends only on $\rho$ up to overall conjugation. 
By definition, $\mathcal{X}(\Gamma)$ is the set of all possible characters.
By invariant theory, $\mathcal{X}(\Gamma)$ has
the structure of an affine complex algebraic variety. Locally, points are determined by finitely many characters from a sufficiently large generating set for $\Gamma$. For more details, we refer to \cite[Cor. 1.4.5]{cullershalen1983varieties}.

The Morgan--Shalen compactification $\overline{\mathcal{X}(\Gamma)}$ can be described as follows \cite{morganshalen1984valuationstrees}:
Let $\mathbb{R}^C$ be the vector space of real-valued functions on $C$ with the weak topology, and $\mathbb{P}(C) := \mathbb{P}(\mathbb{R}^{C})$ the projective space of nonzero real-valued functions on $C$ with the quotient topology. Define a map $\vartheta: \mathcal{X}(\Gamma) \to \mathbb{R}^C$ by
$$
	\vartheta(\rho) :  C  \to   \mathbb{R} \ , \ 
	 [\gamma] \mapsto   \log (|\chi_{\rho}(\gamma)| + 2).
$$
Since $\vartheta(\rho)(1)\neq 0$, the element $\vartheta(\rho)$ has an image in $\mathbb{P}(C)$,  which we denote by $[\vartheta(\rho)]$.
Let $\mathcal{X}(\Gamma)^+$ denote the one-point compactification of $\mathcal{X}(\Gamma)$ with the inclusion map $\iota: \mathcal{X}(\Gamma) \to \mathcal{X}(\Gamma)^+$. The Morgan--Shalen compactification $\overline{\mathcal{X}(\Gamma)}$ is then the closure of the embedded image of $\mathcal{X}(\Gamma)$ in $\mathcal{X}(\Gamma)^+ \times \mathbb{P}(C)$ by the map $\iota \times [\vartheta]$. 

Note that the group $\mathrm{PSL}_2(\mathbb{C}) := \SLC / \mathbb{Z}_2$ is
equal to the identity component of the
isometry group of the 3--dimensional hyperbolic space form $\mathbb{H}^3$. For each representation $\rho$, we define
$$
	\ell_{\rho} :  C  \to   \mathbb{R} \ , \ 
	 [\gamma] \mapsto   \inf_{x \in \mathbb{H}^3} d_{\mathbb{H}^3}(x, \rho(\gamma)x).
$$
It is straightforward to verify  that $|\ell_\rho(\gamma) - 2 \log |\chi_{\rho}(\gamma)|| \leq 2$ for all $\gamma$ (cf.\ \cite{cullershalen1983varieties}). Hence, if a sequence $\rho_n\in \mathcal{X}(\Gamma)$ converges to a point $\rho_\infty$ in $\overline{\mathcal{X}(\Gamma)}\setminus \mathcal{X}(\Gamma)$, then the limit of $[\vartheta(\rho_n)]$ is the same as the limit of $[\ell_{\rho_n}]$. 

We record the following result for later reference.

\begin{lemma}[cf. \cite{cullershalen1983varieties, otal15surveycompactification}]
	\label{lem_MS_cpt_Haus}
$\overline{\mathcal{X}(\Gamma)}$ is Hausdorff and compact.
\end{lemma}

\begin{proof}
	 Since both  $\cX(\Gamma)^+$ and  $\mathbb{P}(C)$  are Hausdorff, the space $\overline{\mathcal{X}(\Gamma)}$ is Hausdorff.
	By \cite[Prop.\ 8]{otal15surveycompactification}, the projection of the closure of the image of $\iota \times [\vartheta]$ to $\mathbb{P}(C)$ is compact. Since $\cX(\Gamma)^+$ is compact, this implies that  $\overline{\mathcal{X}(\Gamma)}$ is compact.\end{proof}

The following result relates the limit of length functions with $\mathbb{R}$--trees.
\begin{theorem}[{\cite{morganshalen1984valuationstrees,morganshalen88degeneration,morganshalen88degernerationtwo}}]
	\label{thm_Morgan-Shalen_compactness}
	Given a sequence of representations $\rho_n: \Gamma \to \mathrm{SL}_2(\mathbb{C})$, the following holds:
	\begin{enumerate}
		\item [(i)] Suppose $\chi_{\rho_n}(\gamma)$ is bounded for each $\gamma \in \Gamma$. Then, after passing to a subsequence if necessary, $\rho_{n}$ converges to an element $\rho_{\infty}$ in $\mathcal{X}(\Gamma)$.
		\item [(ii)] Suppose $\lim_{n \to \infty} |\chi_{\rho_{n}}(\gamma)| = \infty$ for some $\gamma \in \Gamma$. Then, after passing to a subsequence, there exists a minimal $\Gamma$--tree $(\mathcal{T}, d_{\mathcal{T}})$, where the action is denoted by $\rho_{\infty}: \Gamma \to \mathrm{Isom}(\mathcal{T})$, such that its length function $\ell_{\rho_{\infty}}$ satisfies $\ell_{\rho_{\infty}} \not\equiv 0$, and $[\ell_{\rho_n}] \to [\ell_{\rho_{\infty}}]$ in $\mathbb{P}(C)$.
	\end{enumerate}
\end{theorem}

Not every element $[\ell] \in \mathbb{P}(C)$ can be realized as a length
function of a tree (see \cite{chiswell76abstractlengthfunctions} and
\cite[p.\ 586]{cullermorgan97group}). If an element of $\mathbb{P}(C)$ is
realized by the length function of a $\Gamma$--tree, then it can be
realized by a minimal $\Gamma$--tree \cite[Prop.\ 3.1]{cullermorgan97group}. Define $\MPL(\Gamma)$ to be the subspace of $\mathbb{P}(C)$ consisting of all elements that can be realized by the length functions of minimal $\Gamma$--trees that are complete as metric spaces. For each $[\ell] \in \MPL(\Gamma)$, since $\ell$ is not identically zero, the $\Gamma$--action on the corresponding tree has no fixed point.

\section{$\ZT$ Harmonic 1-forms and Taubes' compactness results}
\label{sec_Z2_forms}

In this section, we review the compactness results on stable $\SLC$ flat
connections due to Taubes
\cite{Taubes20133manifoldcompactness,taubes2013compactness}. We will use
these results to define a compactification of the moduli space of flat $\SLC$ connections using $\mathbb{Z}/2$ harmonic 1-forms.

\subsection{$\ZT$ harmonic 1-forms.}
\label{subsec_ZT_1_form}
Let $(M,g)$ be a Riemannian manifold with metric $g$. Here, $M$ is allowed to be non-compact or non-complete
but we assume that $M$ has no boundary. The concept of a
\emph{$\mathbb{Z}/2$ harmonic $1$--form} was first introduced by Taubes
\cite{Taubes20133manifoldcompactness,taubes2013compactness} in order to
describe a compactification of the moduli spaces of solutions to a class of gauge-theoretic equations.

\begin{definition}
	\label{def_ZTharmonicform}
	A $\mathbb{Z}/2$ harmonic $1$--form on $M$ is given by a closed subset $Z \subsetneq M$ and a  two-valued section $\bv$ of $T^*M$ on the complement of $Z$, such that the following conditions hold:
	\begin{itemize}
		\item [(i)] For each $x \notin Z$, there exists an open neighborhood $U \subset M \setminus Z$ of $x$ such that the values of $\bv$ on $U$ have the form $\pm v$, where $v$ is a non-vanishing section of $T^*M|_U$ such that $dv = 0$, $d^*v = 0$.
		\item [(ii)] For every open subset $U$ of $M$ such that $\overline{U}$ is compact, we have $\int_{U \setminus Z} |\bv|^2 < +\infty$, and $\int_{U \setminus Z} |\nabla \bv|^2 < +\infty$.
		\item [(iii)] Let $n$ be the dimension of $M$. There exist constants $C$ and $\epsilon > 0$ such that for every $x \in Z$ and every $r<1$ less than the injectivity radius of $M$ at $x$, we have $\int_{B_r(x)} |\bv|^2 < C \cdot r^{n+\epsilon}$.
	\end{itemize}
\end{definition}

Here, $|\nabla \bv|$ and $|\bv|$ are defined pointwise as the values of $|\nabla v|$ and $|v|$ where $v$ satisfies Condition (i). Condition (iii) above is needed for the proofs of some key properties of $\mathbb{Z}/2$ harmonic $1$--forms. The set $Z$ is called the \emph{zero locus} of $\bv$.

\begin{example}[{\cite[p.\ 9]{Taubes20133manifoldcompactness}}]
	\label{ex_quadraticdifferential}
	When $M$ is a Riemann surface, a $\mathbb{Z}/2$ harmonic 1-form is given by the real part of the square root of a
    holomorphic quadratic differential.
    More explicitly, given a quadratic differential $q \in H^0(K_M^2)$, where $K_M$ is the canonical bundle, $\bv := \Re(\sqrt{q})$ defines a $\mathbb{Z}/2$ harmonic 1-form. Conversely, given a $\mathbb{Z}/2$ harmonic 1-form, let $\bv^{(1,0)}$ be the $(1,0)$ component of $\bv$, then $\bv^{(1,0)} \otimes \bv^{(1,0)}$ also defines a quadratic differential. This correspondence is a bijection. In particular, there exists no $\mathbb{Z}/2$ harmonic 1-form on $S^2$.
\end{example}

We have the following analytic property for the zero loci of $\ZT$ harmonic 1-forms. 
\begin{theorem}[{\cite{taubes2014zero,zhang2017rectifiability}}]
	\label{thm_rectifiability}
	Assume $M$ is connected and the dimension of $M$ is $n\le 4$. Then the zero locus of every $\mathbb{Z}/2$ harmonic $1$--form on $M$ is $(n-2)$--rectifiable with locally finite $(n-2)$--Hausdorff measure.
\end{theorem}

In particular, the zero locus $Z$ of every $\mathbb{Z}/2$ harmonic $1$--form has Lebesgue measure zero. Therefore, we will sometimes write $\int_{U\setminus Z}$ as $\int_U$ when the integrand is defined on the complement of $Z$ and extends to an $L^1$ function over $U\subset M$.

It is convenient to regard a $\mathbb{Z}/2$ harmonic 1-form as a section of the bundle $T^*M/\{\pm 1\}$ over $M$, where the fiber of $T^*M/\{\pm 1\}$ at each point $x \in M$ is the quotient space $(T^*M|_x) / \{\pm 1\}$. A distance function on the quotient space $(T^*M|_x) / \{\pm 1\}$ is defined by 
\begin{equation}
	\label{eqn_defn_2_valued_distance}
	d(\pm v_1, \pm v_2) = \min \{|v_1 - v_2|, |v_1 + v_2|\}.
\end{equation}
We define the moduli space of $\ZT$ harmonic 1-forms as
\begin{equation}
	\label{eqn_MMZ_def}
	\MMZ = \{\bv \mid \bv \; \text{as in Definition \ref{def_ZTharmonicform}, and } \|\bv\|_{L^2(M)} = 1\},
\end{equation}
and endow $\MMZ$ with the $L^2$ topology induced from the distance function \eqref{eqn_defn_2_valued_distance}. 

\begin{proposition}
	\label{prop_MMZ_compact}
	Assume $M$ is closed and has dimension no greater than 4.  Then the space $\MMZ$ is compact.
\end{proposition}

\begin{proof}
	We first show that there exists a constant $C$ depending only on $M$ and the Riemannian metric, such that for each $\bv\in \MMZ$ with zero locus $Z$, we have $\int _{M\setminus Z} |\nabla \bv|^2\le C.$
	
	By the Weitzenb\"ock formula, near every $x\in M\setminus Z$, if we write $\bv$ as $\{\pm v\}$ such that $dv = 0, d^*v=0$, then
	\begin{equation}
		\label{eqn_Weitzenbock_harmonic_form}
		\frac12 d^*d|v|^2 + |\nabla v|^ 2 + \Ric(v,v) = 0,
	\end{equation}
	where $\Ric$ is the Ricci curvature of $M$.  
	
	We show that 
	\begin{equation}
		\label{eqn_int_laplacian_|v|^2}
		\int_{M\setminus Z} d^*d|\bv|^2 =0.
	\end{equation}	
	For each positive integer $i$, let $\rho_i$ be a smooth function on $\mathbb{R}$ such that $\rho_i(x)=0$ when $x\le -i$, $\rho_i(x)=1$ when $x\ge -i+2$, and $\rho_i'(x)\in [0,1]$ for all $x$. Define $\eta_i$ by 
	$$
	\eta_i = \begin{cases}
		\rho_i(\ln |\bv|) & \text{ if } |\bv|\neq 0\\
		0 & \text{ if } |\bv| = 0,
	\end{cases}
	$$
	then $\eta_i$ is a smooth function on $M$. 
	
	On $M\setminus Z$, we have
	$$
	|\nabla \eta_i|  = |\rho_i'(\ln |\bv|)| \cdot \frac{\nabla |\bv|}{|\bv|} \le  \frac{\nabla |\bv|}{|\bv|} \le \frac{|\nabla \bv|}{|\bv|}.
	$$
	Hence the following inequality holds:
	$$
	|\nabla\eta_i|\cdot |\bv|\le |\nabla\bv|.
	$$
	Also note that the support of $\nabla \eta_i$ is a subset of $\{x\in M:|\bv(x)|\in [e^{-i}, e^{-i+2}]\}$. Since $|\nabla \bv|\in L^2(M\setminus Z)$, we conclude that 
	$$
	\lim_{i\to\infty} \int_{\textrm{supp } \eta_i} |\nabla \bv|^2 = 0,
	$$
	and hence $|\nabla \eta_i|\cdot |\bv|$ converges to zero in $L^2(M\setminus Z)$ as $i$ goes to infinity. 
	
	By integration by parts, we have 
	\begin{equation}
		\label{eqn_int_by_parts_d^*d|v|^2}
		\int_{M\setminus Z} (d^*d|\bv|^2 )\cdot \eta_i = \int_{M\setminus Z} (d|\bv|^2 ) (d\eta_i).
	\end{equation}
	By \eqref{eqn_Weitzenbock_harmonic_form}, the function $d^*d|\bv|^2$ is integrable on $M\setminus Z$, and hence the limit of the left-hand side of \eqref{eqn_int_by_parts_d^*d|v|^2} equals $\int_{M\setminus Z} d^*d|\bv|^2$. Since $\|\nabla \bv\|_{L^2(M\setminus Z)} < +\infty $ and $|\nabla \eta_i|\cdot |\bv|\to 0 $ in $L^2$ as $i$ goes to infinity, the right-hand side of \eqref{eqn_int_by_parts_d^*d|v|^2} converges to zero as $i$ goes to infinity. Hence \eqref{eqn_int_laplacian_|v|^2} is proved. By \eqref{eqn_Weitzenbock_harmonic_form}, we conclude that $\int _{M\setminus Z} |\nabla \bv|^2\le C$ for some constant $C$ depending only on $M$.
	
	The Sobolev space $W^{1,p}$ for multi-valued sections of a vector
    bundle was studied in \cite[Ch.\ 4]{de2011q}.  By \cite[Lem.\
    2.1]{zhang2017rectifiability}, the $\mathbb{Z}/2$ harmonic 1-form $\bv$
    can be regarded as a $2$--valued section with $W^{1,2}$ regularity over
    $M$, and its $W^{1,2}$--norm is equal to (up to constant multiplicative
    factors) $\int_{M\setminus Z} |\bv|^2 + |\nabla \bv|^2$.  The desired
    proposition then follows from the Rellich compactness theorem for
    multi-valued sections \cite[Prop.\ 4.6(i)]{de2011q}.
\end{proof}

By \cite[Lem.\ 4.6]{taubes2014zero}, every $\bv\in \MMZ$ is H\"older continuous.

\begin{proposition}
	Assume $M$ is closed. Then the $L^2$ topology and the $\MC^0$ topology coincide on $\MMZ$.
\end{proposition}

\begin{proof}
	Since both topologies are metric spaces, we only need to show that they define the same convergence condition for sequences. It is obvious that $\MC^0$ convergence implies $L^2$ convergence; we show that the converse also holds. Namely, assuming $\{\bv_i\}$ is a sequence in $\MMZ$ that converges to $\bv$ in $L^2$, we show that $\{\bv_i\}$ converges to $\bv$ in $\MC^0$. 
	
	Let $\inj(M)$ denote the injectivity radius of $M$. By \cite[Lem.\ 2.3]{taubes2014zero} and a rescaling argument, there exists a constant $C\ge 1$ depending only on $M$ such that 
	\begin{equation}
		\label{eqn_C0_estimate_from_L2_MMZ}
		\sup_{B_{r/2}(x)} |\bu|^2 \le C\fint_{B_r(x)} |\bu|^2
	\end{equation}
	for all $\bu \in \MMZ$, $x \in M$, and $r < \inj(M)$. 
	
	Let $\{\bv_i\}$ be a sequence in $\MMZ$ that converges to $\bv$ in the $L^2$ topology. 
	Let $D>1$ be a constant such that $\vol \big(B_d(x)\big)< D\cdot \vol \big( B_{d/2}(x) \big)$ for every $x\in M$ and $d<\inj(M)$. 
	For each $\epsilon > 0$, let $B_\epsilon$ be the set consisting of all $x \in M$ such that $|\bv(x)|^2 < \epsilon^2/(8CD)$, where $C$ is the constant in \eqref{eqn_C0_estimate_from_L2_MMZ}. Let $d_\epsilon$ be the distance between $Z$ and $M \setminus B_\epsilon$. Take $\epsilon$ sufficiently small so that $d_\epsilon <\inj (M)$. 
	
	Since $\{\bv_i\}$ converges to $\bv$ in $L^2$,  for $i$ sufficiently large, we have 
	$$
	\fint_{B_{d_\epsilon}(x)}|\bv_i|^2 < \frac{\epsilon^2}{4CD}
	$$
	for all $x \in Z$. Let $V_\epsilon = \cup_{x \in Z} B_{d_\epsilon/2}(x)$. By \eqref{eqn_C0_estimate_from_L2_MMZ}, this implies
	$$
	\sup_{V_\epsilon}|\bv_i|^2 < \frac{\epsilon^2}{4}.
	$$
	 So the $\MC^0$--distance between $\bv_i$ and $\bv$ on $V_\epsilon$ is no greater than 
	$$
	\sup_{V_\epsilon} (|\bv_i| + |\bv|) < \epsilon (1/2 + \sqrt{1/(8CD)}) < \epsilon. 
	$$
	On the other hand, $\bv_i$ converges to $\bv$ uniformly on $M \setminus V_\epsilon$ by standard elliptic estimates. Hence the desired result is proved.
\end{proof}

\subsection{The moduli space $\MM_{\SLC}$ and its compactification}
\label{subsec_BMMc}

Let $M$ be a closed oriented 3-manifold. Let $P$ be a principal $\SU(2)$ bundle over $M$, and let $\mathfrak{g}_P$ be the associated $\mathfrak{su}(2)$ bundle given by the adjoint action. Consider the system of equations 
\begin{equation}
	\label{eq_flat_connection}
	\begin{split}
		F_A = \frac{1}{2} [\phi, \phi], \quad d_A\phi=0, \quad d_A^*\phi=0,
	\end{split}
\end{equation}
where $A$ is a connection on $P$ and $\phi$ is a section of $T^*M \otimes \mathfrak{g}_P$. We define $\MMc$ to be the set of solutions to \eqref{eq_flat_connection}, 
modulo $\SU(2)$ gauge transformations. The topology on $\MMc$ is given by the $W^{k,p}$--Sobolev norm for $k, p$ sufficiently large. By the standard elliptic bootstrapping argument, the topology on $\MMc$ does not depend on the choice of $(k,p)$ when $p > 1$ and $k$ is sufficiently large. 

Note that if $\bv$ is a $\mathbb{Z}/2$ harmonic $1$--form, then $\bv \otimes \bv$ is a single-valued section of $T^*M \otimes T^*M$. 
The following result is a consequence of Taubes' compactness theorems.

\begin{theorem}[\cite{Taubes20133manifoldcompactness,walpuskizhang2019compactness,parker2023concentrating}]
	\label{thm_taubes_compactness}
	Let $(A_i, \phi_i)$ be a sequence of solutions to \eqref{eq_flat_connection}, and let $r_i := \|\phi_i\|_{L^2(M)}$.
	\begin{itemize}
		\item [(i)] If $\{r_i\}$ is bounded, then there exists a subsequence of $(A_i, \phi_i)$ that converges in $\MC^\infty$ after gauge transformations.
		\item [(ii)] If $\lim_{i\to\infty} r_i = +\infty$, then there exists a subsequence (which we still denote by $(A_i, \phi_i, r_i)$) and $\bv \in \MMZ$ such that $r_i^{-2}\Tr(\phi_i \otimes \phi_i)$ converges to $\bv \otimes \bv$ in $\MC^0$.
	\end{itemize}
\end{theorem}

Theorem \ref{thm_taubes_compactness} is implicitly contained in \cite{Taubes20133manifoldcompactness, walpuskizhang2019compactness, parker2023concentrating}. In the following, we deduce the statement of Theorem \ref{thm_taubes_compactness} from the above references. 
\begin{proof}
	[Proof of Theorem \ref{thm_taubes_compactness}]
	Case (i) follows from standard elliptic bootstrapping. For case (ii), results in \cite{Taubes20133manifoldcompactness, walpuskizhang2019compactness, parker2023concentrating} imply that there exists a subsequence of $(A_i, \phi_i, r_i)$, which we denote by the same notation, and a $\bv \in \MMZ$ with zero locus $Z$, such that 
	\begin{enumerate}
		\item $|\phi_i|/r_i$ converges to $|\bv|$ in $\MC^0$.
		\item For every compact set $K$ contained in an open ball in $M \setminus Z$, there exists $\phi$ on $K$ such that, after a sequence of gauge transformations, $\phi_i/r_i$ converges to $\phi$ in the weak $W^{2,2}$ topology on $K$.
		\item The spinor $\phi$ in (ii) satisfies $\Tr(\phi \otimes \phi) = \bv \otimes \bv$. 
	\end{enumerate}
	Statements (1) and (2) above are directly given by \cite[Thm.\
    1.28]{walpuskizhang2019compactness}. A stronger result was proved in
    \cite[Thm.\ 1.3]{parker2023concentrating}, where $\phi_i/r_i$ was
    shown to converge to $\phi$ in the $\MC^\infty$ topology on $K$. A
    similar compactness result was given in \cite[Thm.\ 1.1a]{Taubes20133manifoldcompactness}
    under weaker assumptions, but with weaker Sobolev regularity in the convergence statement.
    Statement (3) above follows from part (3) of the second bullet point of
    \cite[Thm.\ 1.1a]{Taubes20133manifoldcompactness}.
	
	Now we prove Case (ii) of Theorem \ref{thm_taubes_compactness}. 
	For each $\epsilon > 0$, let $B_\epsilon$ be the set of $x \in M$ such that $|\bv \otimes \bv| \leq \epsilon/4$. Then $B_\epsilon$ is a closed subset of $M$ that contains an open neighborhood of $Z$. By Statement (1) above, for $n$ sufficiently large, we have
	$$\sup_{B_\epsilon} r_i^{-2}|\Tr (\phi_i \otimes \phi_i)| < \epsilon/2.$$
	Hence the $\MC^0$ distance between $r_i^{-2}|\Tr (\phi_i \otimes \phi_i)|$ and $\bv \otimes \bv$ is less than $\epsilon$ on $B_\epsilon$ when $i$ is sufficiently large. 
	By Statement (2) above and the Sobolev embedding theorems, after passing to a subsequence if necessary, the sequence $\phi_i/r_i$ converges to $\phi$ on $K$ in the $\MC^0$ topology. Since $\Tr(\phi_i\otimes \phi_i)$ is gauge invariant, the result is proved. 
\end{proof}

Next, we define a compactification of $\MMc$ using Theorem \ref{thm_taubes_compactness}. Let $\mbM$ be the disjoint union of $\MMc$ and $\MMZ$. We define a topology on $\mbM$ as follows. The open sets on $\mbM$ are generated by the following two collections of subsets:
\begin{enumerate}
	\item[(i)] Open subsets of $\MMc$.
	\item[(ii)] For $\bv \in \MMZ$, $\epsilon > 0$, $N > 0$, the subset of $\mbM$ that contains all $\bv'$ such that the $\MC^0$ distance between $\bv \otimes \bv$ and $\bv' \otimes \bv'$ is less than $\epsilon$, and all $(A, \phi) \in \MMc$ such that 
	\begin{enumerate}
		\item $\|\phi\|_{L^2} > N$,
		\item the $\MC^0$ distance between $\|\phi\|_{L^2}^{-2}\Tr(\phi \otimes \phi)$ and $\bv \otimes \bv$ is less than $\epsilon$.
	\end{enumerate}
\end{enumerate}

Note that for $\bv, \bv' \in \MMZ$, we have $\bv = \bv'$ if and only if $\bv \otimes \bv = \bv' \otimes \bv'$, and that the $\MC^0$ topology on $\MMZ$ is the same as the pull-back topology from $\MC^0(T^*M \otimes T^*M)$ via the map 
	$\bv \mapsto \bv \otimes \bv$.
Therefore,  $\MMc$ and $\MMZ$ are homeomorphic to their embedded images in $\mbM$. Moreover, Proposition \ref{prop_MMZ_compact} and Theorem \ref{thm_taubes_compactness} can be summarized into the following statement.

\begin{corollary}
	\label{cor_mbM_cpt_Haus}
	The space $\mbM$ is Hausdorff and compact. 
\end{corollary} 

\begin{proof}
	Since both $\MMc$ and $\MMZ$ are Hausdorff, the fact that $\mbM$ is Hausdorff follows straightforwardly from the definition of its topology. Proposition \ref{prop_MMZ_compact} and Theorem \ref{thm_taubes_compactness} imply that $\mbM$ is sequentially compact. It is also straightforward to verify that $\mbM$ is first countable. Hence $\mbM$ is compact. 
\end{proof}

\begin{definition}
	\label{def_BMMc}
	Define $\BMMc$ to be the closure of $\MMc$ in $\mbM$. 
\end{definition}
We call $\BMMc$ the \emph{compactification} of $\MMc$. 
\begin{remark}
	One can construct a more refined compactification of $\MMc$ from the analytical results in \cite{Taubes20133manifoldcompactness, haydyswalpuski2015compactness, walpuskizhang2019compactness, parker2023concentrating} by also considering the limits of the connection terms $A_i$. However, it is not clear to us how the convergence of the connection terms is related to the Morgan--Shalen compactification. See \cite{ott2020higgs,he2023algebraic} for results in this direction in the two-dimensional case. 
\end{remark}

\section{Measured foliations and $\ZT$ harmonic 1-forms}
\label{sec_leaf_space}

In this section, we review a construction of measured foliations from
$\mathbb{Z}/2$ harmonic $1$--forms by Taubes
 \cite[p.\ 14]{Taubes20133manifoldcompactness}. 
 Then we prove that on simply connected manifolds, the leaf space of the measured foliation defined from a $\mathbb{Z}/2$ harmonic $1$--form is always an $\mathbb{R}$--tree.

\subsection{Measured foliations defined by $\mathbb{Z}/2$ harmonic 1-forms}
\label{subsec_ZT_form_to_foliation}
The theory of measured foliations and quadratic differentials on Riemann
surfaces has found significant applications in the geometry and topology of
Riemann surfaces and 3-manifolds 
\cite{Hubbard1, Hubbard2, thurston97threedimensionlgeometry}.
For three- or four-dimensional Riemannian manifolds, concepts
such as singular measured foliations, measured laminations, and weighted
branched surfaces have been extensively developed in
\cite{hatcher1996fulllamnination,gabai89essential,oertel88measured}. We
review the construction of a singular measured foliation from a
$\mathbb{Z}/2$ harmonic $1$--form, following Taubes \cite[p.\ 14]{Taubes20133manifoldcompactness}. 

We first introduce the concept of a singular measured foliation.
\begin{definition} \label{def:singular-foliation}
	Let $M$ be a manifold with dimension $n$.
	Let $Z$ be a closed subset of Hausdorff codimension at least $2$. A (codimension-one) \emph{singular foliation} $\mathcal{F}$ on $M$ with singular set $Z$ is a smooth foliation on $M \setminus Z$ with codimension-one leaves. A \emph{transverse measure} $\mu$ on $\mathcal{F}$ is  a measure for arcs in $M\setminus Z$ such that the following conditions hold:
	\begin{enumerate}
		\item [(i)] $\mu$ is non-zero on transverse arcs.
		\item [(ii)] $\mu$ vanishes on an arc if and only if the arc is tangent to a leaf.
		\item [(iii)] (Holonomy invariance) $\mu$ is invariant along
            homotopies among transverse arcs
		that keep endpoints in the same leaves.
	\end{enumerate}
\end{definition}
We will sometimes refer to a ``singular foliation'' as a ``foliation'' when there is no risk of confusion.

Now we associate a measured foliation to a $\mathbb{Z}/2$ harmonic 1-form.
To provide some intuition, recall from Example
\ref{ex_quadraticdifferential} that, on a Riemann surface, a $\mathbb{Z}/2$
harmonic 1-form is the real part of the square root of a quadratic
differential. Let us briefly review how a quadratic differential gives rise to a measured foliation.

\begin{example}
	Suppose $M$ is a closed Riemann surface. Given a nonzero holomorphic  quadratic differential $q$,
    the zeros $Z := q^{-1}(0)$ form the singular set of a foliation. Away
    from $Z$, one can write $q$ locally as $\omega^2$ for some holomorphic
    $1$--form $\omega$. The kernel of $\textrm{Re}\,\omega$ then defines a
    foliation, with a transverse measure given by  $|\textrm{Re}\,\omega|$. This local construction can be combined to form a global measured foliation, often referred to as the \emph{vertical measured foliation} of $q$.
\end{example}

This construction can be generalized to $\mathbb{Z}/2$ harmonic 1-forms as in \cite{taubes2013compactness}. Consider a Riemannian manifold $M$ (not necessarily closed or complete), and let $\bv$ be a $\mathbb{Z}/2$ harmonic 1-form on $M$. We define the singular set as $Z = |\bv|^{-1}(0)$. By Theorem \ref{thm_rectifiability}, if the dimension of $M$ is no greater than $4$, then $Z$ meets the conditions described in Definition \ref{def:singular-foliation}. 

For a point $x \in M \setminus Z$, let $U_x$ be a small neighborhood around $x$ where $\bv$ can be locally expressed as $\pm v$, with $v$ a single-valued 1-form on $U_x$. By Definition \ref{def_ZTharmonicform}, $v$ is a closed 1-form. The kernel of $v$ thus defines a codimension-one foliation on $U_x$, which is independent of the choice of sign for $v$. Consequently, $\ker \bv$ defines a smooth foliation on $M \setminus Z$.

We define a transverse measure $\mu_\bv$ for this foliation as follows. For each $\MC^1$ arc $\gamma: [0,1] \to M \setminus Z$, the transverse length of $\gamma$ is given by
\begin{equation}
	\label{eqn_transverse_measure_from_bv}
	\mu_\bv(\gamma) = \int_0^1 |\bv(\dot{\gamma}(t))| dt.
\end{equation}
It’s important to note that, similar to the cases of $|\bv|$ and $|\nabla\bv|$, the value of $|\bv(\dot{\gamma}(t))|$ is well-defined despite $\bv$ being a two-valued section. Since $\bv$ is locally given by a closed smooth form in $M \setminus Z$, the measure $\mu_\bv$ is holonomy invariant. As a result, the pair $(\ker \bv, \mu_\bv)$ defines a measured foliation on $M \setminus Z$.

\subsection{The leaf spaces of $\ZT$ harmonic 1-forms}
Next, we define the leaf space associated to a $\ZT$ harmonic form. We
prove that if the background manifold is simply connected, then the leaf space is an $\mathbb{R}$--tree. 

Consider a connected Riemannian manifold $M$ (not necessarily closed or complete), and let $\bv$ be a $\mathbb{Z}/2$ harmonic $1$--form on $M$. We define a pseudo-metric on $M$ by
\begin{equation}
	\label{eq_lengthfunction_defined_by_Z2harmonicform}
	d_{M,\bv}(x,y) = \inf_\gamma \int_0^1 |\bv(\dot{\gamma}(t))| dt,
\end{equation}
where $\gamma$ ranges over all piecewise $\MC^1$ curves from $x$ to $y$.
Note that the definition of $d_{M,\bv}$ depends globally on both $M$ and $\bv$: if $U$ is an open subset of $M$, then it is not necessarily true that $d_{U,\bv}$ equals the restriction of $d_{M,\bv}$ on $U$.

\begin{definition}
	\label{def_leaf_space}
	The metric space $\MT_{M,\bv}$ is the quotient space of $M$ with respect to the pseudo-metric $d_{M,\bv}$.
\end{definition}

We call the space $\MT_{M,\bv}$ the \emph{leaf space} of $\ker \bv$ on $M$, and we will use $d_{M,\bv}$ to denote the metric on this quotient space as well. 

\begin{theorem} \label{thm_pseudo-metric-R-tree}
	Suppose $M$ is a simply connected manifold, $\bv$ is a $\ZT$ harmonic $1$--form on $M$ with zero locus $Z$. If the dimension of $M$ is greater than $4$, we assume in addition that $\bv$ is continuous and $Z$ is rectifiable with locally finite Hausdorff measure in codimension $2$. Then the metric space $(\mathcal{T}_{M,\bv}, d_{M,\bv})$ is an $\mathbb{R}$--tree. 
\end{theorem}

The remainder of this section is devoted to the proof of Theorem
\ref{thm_pseudo-metric-R-tree}.
We first recall the following standard result about $\mathbb{R}$--trees.

\begin{proposition}[{\cite[Prop.\ B.31]{py2023lectures}}]
	A metric space $(X,d)$ is an $\mathbb{R}$--tree if and only if it is path connected and the following inequality holds for all $x,y,z,t\in X$:
	$$
	d(x,z)+d(y,t) \le \max\{d(x,y)+d(z,t),d(x,t)+d(y,z)\}.
	$$
\end{proposition}

Now we prove the following technical lemma. 
\begin{lemma}
	\label{lem_oppo_pts}
	Suppose $D$ is the closed unit disk in $\mathbb{R}^2$, and let $int(D)$ denote the interior of $D$. Suppose $\bv$ is a continuous, two--valued $1$--form defined on an open neighborhood of $D$ in $\mathbb{R}^2$ with zero locus $Z$, such that away from its zero points $\bv$ is locally given by $\pm v$ with $v$ being closed. 
	We further assume that $Z\cap \partial D$ is finite, $Z\cap int(D)$ is compact, and that for each $\delta>0$, there exists a smooth domain $\Omega_\delta\subset int(D)$ such that $\Omega_\delta$ contains $Z\cap int(D)$, and  $\Omega_\delta$ is contained in the $\delta$--neighborhood of  $Z\cap int(D)$, and the total length (i.e. the Hausdorff measure in dimension $1$) of $\partial \Omega_\delta$ is less than $1$.
Finally, 	
	assume that $\ker \bv$ is transverse to $\partial D$ at all but finitely many points.
	
	Let $d_{D,\bv}$ be the pseudo-distance function on $D$ associated with $\bv$ given by \eqref{eq_lengthfunction_defined_by_Z2harmonicform}.  Let $p_1,p_2,p_3,p_4$ be four distinct points appearing in cyclic order on $\partial D$, the points $p_1,p_2,p_3,p_4$ divide $\partial D$ into four arcs. 
	
	Then there exist two points $a$ and $b$ on a pair of opposite (closed) arcs divided by $p_1,p_2,p_3,p_4$ such that $d_{D,\bv}(a,b)=0$.
\end{lemma}

\begin{proof}
	We adapt an argument from \cite[Lem.\ III.4]{levitt93constructing}. To simplify notation, we will write $d_{D,\bv}$ as $d$ in the proof. 
	Assume $p_1,p_2,p_3,p_4$ are ordered counterclockwise.  In the following, for $x,y\in \partial D$, we use $\overline{xy}$ to denote the closed arc in $\partial D$ bounded by $x, y$, where the arc goes from $x$ to $y$ in the counterclockwise direction. If $A,B$ are compact subsets of $D$, we use $d(A,B)$ to denote $\inf_{a\in A,b\in B} d(a,b)$. 
	
	If $d(\overline{p_1p_2},\overline{p_3p_4})=0$, then the desired result already holds. In the following, assume  $d(\overline{p_1p_2},\overline{p_3p_4})>0$. Let $q\in \overline{p_1p_2}$ be such that $d(q,\overline{p_4p_1})=0$ and that $q$ is  furthest away from $p_1$ as a point on the arc $\overline{p_1p_2}$ among all points satisfying this condition. Similarly, let $r\in \overline{p_1p_2}$ be such that $d(r,\overline{p_2p_3})=0$ and that $r$ is furthest away from $p_2$ among all points in $\overline{p_1p_2}$ satisfying this condition.

	We discuss three cases.  If $q=r$, then the above conditions imply that $d(\overline{p_4p_1},\overline{p_2p_3})=0$, and the desired result holds.
	
	If $r\neq q$ and $\overline{rq}\subset \overline{p_1p_2}$ (that is, $q$ is on the counterclockwise side of $r$), then by the above conditions, for every $\epsilon>0$, there exists a $\MC^1$ curve from $r$ to $\overline{p_2p_3}$, and a $\MC^1$ curve from $q$ to $\overline{p_4p_1}$, whose lengths with respect to $d$ are less than $\epsilon$.  These two curves must intersect, so $d(q,r)<2\epsilon$. Since this statement holds for all $\epsilon>0$, we conclude that $d(q,r)=0$, which implies $d(\overline{p_4p_1},\overline{p_2p_3})=0$, and the desired result holds.
	
	If $r\neq q$ and  $\overline{qr}\subset \overline{p_1p_2}$ (that is, $r$ is on the counterclockwise side of $q$), we consider the singular foliation $\mathcal{F}$ defined by $\ker \bv$ on $D$. Since $\bv$ is locally given by closed forms away from zero points, there is a transverse invariant measure $\mu_\bv$ on $\mathcal{F}$ as given by \eqref{eqn_transverse_measure_from_bv}. 
	
	For each point $x\in \overline{qr}\setminus\{q,r\}$, if $\mathcal{F}$ is defined and is transverse to $\partial D$ at $x$, we consider the leaf of $\mathcal{F}$ passing through $x$. View the leaf as a parametrized curve starting at $x$. By the Poincar\'e--Bendixson theorem, one of the following holds:
	\begin{enumerate}
		\item The leaf intersects a boundary point of $\partial D$ other than $x$.
		\item The leaf converges to a zero point of $\bv$.
		\item The leaf converges to a limit cycle.
	\end{enumerate}
	Since $\mathcal{F}$ admits a transverse measure, Case (3) cannot happen. By the definitions of $q,r$, we know that $d(\{x\},\overline{rp_3})>0$, $d(\{x\},\overline{p_4q})>0$. By the assumption that $d(\overline{p_1p_2},\overline{p_3p_4})>0$, we know that $d(x,\overline{p_3p_4})>0$.  As a result, if Case (1) happens, then the intersection point of the leaf with $\partial D$ is in the interior of $\overline{qr}$.
	
	Let $U\subset \overline{qr}\setminus\{q,r\}$ be the set of points $x$ such that 
	\begin{enumerate}
		\item $\bv$ is non-zero at $x$ and $\mathcal{F} = \ker \bv$ is transverse to $\partial D$ at $x$.
		\item The leaf of $\mathcal{F}$ starting at $x$ transversely intersects a point of $\overline{qr}\setminus\{q,r\}$ other than $x$.
	\end{enumerate}
	
	Then $U$ is an open subset of $\overline{qr}\setminus\{q,r\}$. We claim that 
	\begin{equation}
		\label{eqn_U_full_measure}
		\mu_\bv(U) = \mu_\bv(\overline{qr}).
	\end{equation}
 This is because the leafs emanating from the zero points of $\bv$ have zero measure with respect to $\mu_\bv$. More precisely, recall that $Z$ denotes the zero set of $\bv$. For each $\epsilon>0$, there exists $\delta$ such that $|\bv|<\epsilon$ on the $\delta$--neighborhood of $Z$. Let $\Omega_\delta$ be the smooth domain given by the assumptions on $Z$, then $\mu_\bv(\partial \Omega_\delta)<\epsilon$. Then the set of points $x\in \overline{qr}$ such that there is leaf of $\mathcal{F}$ intersecting $\partial D$ transversely at $x$ and enters $\Omega_\delta$ has $\mu_\bv$--measure less than $\epsilon$.
 Similarly, if $p$ is a zero point of $\bv$ in $\partial D$, then $\mu_\bv(\partial B_r(p)\cap D)\to 0$ as $r\to 0$, so the  set of points $x\in \overline{qr}$ such that there is leaf of $\mathcal{F}$ intersecting $\partial D$ transversely at $x$ and converges to $p$ on the other end has $\mu_\bv$--measure zero. Recall that for every point $x\in\overline{qr}\setminus(\{q,r\}\cup U)$, one of the following conditions holds:
 \begin{enumerate}
 	\item $x$ is a zero point of $\bv$,
 	\item $x$ is a point where $\mathcal{F}$ is tangent to $\partial D$, 
 	\item $x$ is a point whose leaf passes through a tangent point of $\mathcal{F}$ with $\partial D$,
 	\item $x$ is a point whose leaf converges to a zero point of $\bv$.
 \end{enumerate}
 The first three cases only contain finitely many points, and points in the last case have measure zero with respect to $\mu_\bv$ by the previous argument. Hence Equation \eqref{eqn_U_full_measure} is proved. Since $\ker\bv$ is transverse to $\partial D$ at all but finitely many points, we conclude that $U$ has full measure in $\overline{qr}$ with respect to the standard Lebesgue measure as well.
	
	Define an involution 
	$
	\iota: U\to U
	$,
	such that for every $x\in U$, the image $\iota(x)\in U$ is the other endpoint of the leaf of $\mathcal{F}$ passing through $x$. The map $\iota$ is orientation-reversing. 
	
	Parameterize the arc $\overline{qr}$ by the interval $[0,1]$ via a smooth diffeomorphism:
	$\varphi:[0,1]\to \overline{qr}$.
	Define a function 
$$
		\xi: [0,1]\to \mathbb{R}: t\mapsto d(q,\varphi(t)),
$$
	where we recall that $d$ denotes the pseudo-distance function $d_{D,\bv}$. 
	Then $\xi$ a Lipschitz function on $[0,1]$, so $\xi'$ exists almost everywhere and 
	$$\xi(1)-\xi(0) = \int_0^1\xi'.$$ 
	For each open interval $(s,t)\in \varphi^{-1}(U)$, let $\hat s = \varphi^{-1}\circ \iota \circ \varphi(s)$ and $\hat t = \varphi^{-1}\circ \iota \circ \varphi(t)$, we have 
	$$
	\int_s^t \xi' = d(q,\varphi(t)) - d(q,\varphi(s)) = d(q,\varphi(\hat t)) - d(q,\varphi(\hat s)) = - \int_{\hat t}^{\hat s} \xi'.
	$$
	Hence the involution $\varphi^{-1}\circ \iota\circ \varphi$ on $\varphi^{-1}(U)$ reverses the orientation and preserves the differential form $\xi'(t)dt$. This implies $\int_{\varphi^{-1}(U)} \xi' =0$. Since $\varphi^{-1}(U)$ has full Lebesgue measure in $[0,1]$, we conclude that $\int_0^1 \xi'=0$,  so $d(q,r)=\xi(1)=0$.  This contradicts the definitions of $q,r$.
\end{proof}

We now prove Theorem \ref{thm_pseudo-metric-R-tree}.

\begin{proof}[Proof of Theorem \ref{thm_pseudo-metric-R-tree}]
	Let $p_1,p_2,p_3,p_4$ be distinct points in $M$. By Lemma \ref{lem_oppo_pts}, we only need to show that 
	\begin{equation}
		\label{eqn_quadriple_ineq}
		\begin{split}
		&d_{M,\bv}(p_1,p_3) + d_{M,\bv}(p_2,p_4) \\
		\le &\max\{d_{M,\bv}(p_1,p_2) + d_{M,\bv}(p_3,p_4) , d_{M,\bv}(p_1,p_4) + d_{M,\bv}(p_2,p_3) \}.
		\end{split}
	\end{equation}
	Since \eqref{eqn_quadriple_ineq} is a closed condition on $(p_1,p_2,p_3,p_4)$, we may assume without loss of generality that $p_1,p_2,p_3,p_4$ are not zero points of $\bv$. For notational convenience, we will interpret the subscripts modulo $4$. Fix $\epsilon >0 $. For each $i$, let $\gamma_i$ be a smooth arc from $p_i$ to $p_{i+1}$ such that 
	$$
	\mu_\bv(\gamma_i)\le d_{M,\bv}(p_i,p_{i+1})+\epsilon. 
	$$
	After perturbation, we may assume that the union of $\gamma_1,\dots,\gamma_4$ forms a smoothly embedded circle in $M$ that is disjoint from the zero locus of $\bv$. Denote this circle by $C$. After a further perturbation if necessary, we may assume that $C$ is transverse to $\ker\bv$ at all but finitely many points. Since $M$ is simply connected, there exists a smooth immersion 
	$
	f:D\to M
	$,
	where $D$ is the unit disk in $\mathbb{R}^2$, such that $f(\partial D)= C$. Since the zero locus of $\bv$ has locally finite codimension two Hausdorff measure, we may perturb $f$, with $f|_{\partial D}$ fixed, so that $f(D)$ intersects the zero locus of $\bv$ at finitely many points, and that the zero locus of $f^*(\bv)$ is discrete on the complement of $(f\circ |\bv|)^{-1}(0)$. As a result, the assumptions of Lemma \ref{lem_oppo_pts} holds for $f^*(\bv)$. 
	
	Applying Lemma \ref{lem_oppo_pts} to $f^{-1}(p_i)$ for $i=1,2,3,4$ and to the $2$--valued form $f^*(\bv)$, we conclude that there exists $a,b\in \partial D$ on opposite arcs such that $d_{D,f^*(\bv)}(a,b)=0$. Since $d_{M,\bv}(a,b)\le d_{D,f^*(\bv)}(a,b)$, this implies $d_{M,\bv}(a,b)=0$. Assume without loss of generality that $a$ is in the arc bounded by $p_1,p_2$, and $b$ is in the arc bounded by $p_3, p_4$. Then we have
	$$
	d_{M,\bv}(p_1,p_3) \le d_{M,\bv}(p_1,a) + d_{M,\bv}(b,p_3),\qquad d_{M,\bv}(p_2,p_4) \le d_{M,\bv}(p_2,a) + d_{M,\bv}(b,p_4).
	$$
	By the assumptions on $f(\partial D)$, we have
	$$
	d_{M,\bv}(p_1,a) + d_{M,\bv}(p_2,a) \le \mu_\bv(\gamma_1) \le  d_{M,\bv}(p_1,p_2) + \epsilon,
	$$
	and similarly,
	$$
	d_{M,\bv}(p_3,b) + d_{M,\bv}(p_4,b) \le  d_{M,\bv}(p_3,p_4) + \epsilon.
	$$
	As a result,
	$$
	d_{M,\bv}(p_1,p_3) + d_{M,\bv}(p_2,p_4) \le  d_{M,\bv}(p_1,p_2) + d_{M,\bv}(p_3,p_4) + 2\epsilon.
	$$
	The desired result then follows by taking $\epsilon \to 0$. 
\end{proof}



\section{Harmonic maps to $\mathbb{R}$--trees}
\label{sec_harmonicmaps}

Let $M$ be a closed $3$--manifold and $\widetilde{M}$  its universal cover.
In this section, we show that if $u:\tM\to \MT$ is a
$\pi_1(M)$--equivariant harmonic map to an $\mathbb{R}$--tree, then the
gradient of $u$ defines a $\mathbb{Z}/2$ harmonic 1-form on $M$. We first
recall the formulation  of Korevaar--Schoen which  gives meaning to  the
gradient of $u$ as a Radon--Nikodym derivative.
Then, we use regularity results of harmonic maps to construct the associated $\mathbb{Z}/2$ harmonic $1$--form.

\subsection{Energy density and directional derivatives}
\label{subsec_energy_density}
Here, we review some  terminology from \cite{korevaar1993sobolev}, specialized to the case of $L^2$ norms. Following the notation of  \cite{korevaar1993sobolev},
let $(\Omega,g)$ be an $n$-dimensional Riemannian manifold and $(X,d)$ a complete metric space. 
Unless otherwise specified, $\Omega$ is allowed to be non-compact or non-complete but we assume that $\Omega$ has no boundary.
Later, we will take $\Omega$ to be the universal cover of a closed $3$--manifold and $(X,d)$ to be an $\mathbb{R}$--tree.

For $\epsilon>0$, let $\Omega_\epsilon$ be the set of all points $x\in \Omega$ such that the exponential map at $x$ is well-defined on the open ball with radius $\epsilon$. Let $\mu$ denote the volume measure of $\Omega$.

A map $u:\Omega\to X$ is called \emph{locally $L^2$} if for every compact subset $K$ in $\Omega$ and every $p\in X$, we have 
$$
\int_K d^2(u(x),p) \,d\mu(x) <+\infty.
$$
If $u$ is locally $L^2$, define the \emph{$\epsilon$--approximate energy function} $e_{\epsilon}: \Omega_\epsilon \to \mathbb{R}$ as
\begin{equation}
	\label{den_e_epsilon}
	e_{\epsilon}(x) := \frac{1}{\omega_{n-1}\epsilon^{n-1}} \int_{\partial B_{\epsilon}(x)} \frac{d^2(u(x), u(y))}{\epsilon^2} \, d\sigma,
\end{equation}
where $\omega_{n-1}$ is the area of $S^{n-1}$ in $\mathbb{R}^n$ (recall that $n$ is the dimension of $\Omega$), and $d\sigma$ is the area form on $\partial B_\epsilon(x)$. The map $u$ is said to have \emph{finite energy}, if 
$$
\sup_{\varphi\in \MC_0^\infty(M),0\le\varphi\le 1} \limsup_{\epsilon\to 0} \int_\Omega e_\epsilon\cdot \varphi \,d\mu<+\infty.
$$
The space of maps from $\Omega$ to $X$ with finite energy is denoted by
$W^{1,2}(\Omega,X)$. By \cite[Thm.\  1.10]{korevaar1993sobolev}, if $u\in W^{1,2}(\Omega,X)$, then the measures $e_{\epsilon}(x) \mu(x)$ converge weakly to a limit $e(x)\mu(x)$ as $\epsilon\to 0$, where $e(x)\in L^1(\Omega)$. The function $e(x)$ is called the \emph{energy density function}, and we define $|\nabla u|(x) := e(x)^{1/2} \in L^2(\Omega)$.

The directional derivatives are defined in a similar way.  Assume $Z$ is a Lipschitz vector field over $\Omega$. For $x\in \Omega$, let $x+\epsilon Z$ denote the image of $x$ after flowing along $Z$ for time $\epsilon$. Define $^{Z}e_{\epsilon}(x) = d^2(u(x),u(x+\epsilon Z))/\epsilon^2$. By \cite{korevaar1993sobolev}, if $u\in W^{1,2}(\Omega,X)$, then there exists a unique non-negative function $|u_*(Z)|\in L^2(\Omega)$ such that the measures $^{Z}e_{\epsilon}(x) \mu(x)$ converge weakly to $|u_*(Z)|^2(x)\mu (x)$.

A map $u\in W^{1,2}_{loc}(\Omega,\mathbb{R})$ is called \emph{harmonic}, if
it is a critical point of the Dirichlet functional (we refer the reader to
\cite[Sec.\ 2.2]{korevaar1993sobolev} for more details). If the target $(X,d)$ is an $\mathbb{R}$--tree, then harmonic maps to $X$ minimizes the energy with respect to compactly supported perturbations. 

Now assume $M$ is a closed manifold and let $\widetilde{M}$ be its universal cover.  Assume $(\MT,d_\MT)$ is a complete $\mathbb{R}$--tree with a fixed isometric action by $\pi_1(M)$. 
The following results established the existence and uniqueness of
equivariant harmonic maps. The existence result follows from the work of
Korevaar--Schoen \cite{korevaar1993sobolev, korevaari1997global}, and the uniqueness follows from \cite{mese2002uniqueness}. 
See also \cite{daskalopoulos2021uniqueness} for generalizations.

%

\begin{theorem}[\cite{korevaar1993sobolev, korevaari1997global,mese2002uniqueness,daskalopoulos2021uniqueness}]
	\label{thm_existence_uniqueness}
	Suppose the $\pi_1(M)$ action on a complete $\mbR$ tree $(\mathcal{T},d_{\MT})$ has no fixed ends (see the definitions above Theorem \ref{thm_nonabelian_determinant_tree}). Then there exists a $\pi_1(M)$--equivariant harmonic map $u: \widetilde{M} \to \mathcal{T}$. Moreover, if $u_0, u_1: \widetilde{M} \to \mathcal{T}$ are two  $\pi_1(M)$--equivariant harmonic maps such that $u_0\neq u_1$, then either $u_0(\widetilde{M})$ or $u_1(\widetilde{M})$ is contained in a geodesic.
\end{theorem}

Mese \cite{mese2002uniqueness} showed  that distinct harmonic maps always have the same directional derivatives. 

\begin{proposition}[{\cite[Cor.\ 13]{mese2002uniqueness}}]
	\label{prop_uniqueness_directional_energy}
	Fix a $\pi_1(M)$ action on an $\mathbb{R}$--tree $(\MT,d_\MT)$. 
	Let $u_0, u_1: \widetilde{M} \to \mathcal{T}$ be two equivariant harmonic maps to $\MT$, and let $Z$ be a Lipschitz vector field on $\widetilde{M}$. Then we have $|(u_0)_* Z| = |(u_1)_* Z|$ (almost everywhere).
\end{proposition}

\subsection{Regularity of harmonic maps to $\mathbb{R}$--trees}

Now we review some regularity results on harmonic maps to $\mathbb{R}$--trees from the literature.

\begin{definition}
	\label{def_regular_set}
	Let $u: \Omega \to \mathcal{T}$ be a harmonic map from a manifold to a tree. 
	We say that $p\in \Omega$ is a \emph{regular point} of $u$, if there exists an open ball $B_r(p)$ centered at $p$ such that $u(B_r(p))$ is contained in a geodesic in $\mathcal{T}$. If $p$ is not a regular point, then it is called a \emph{singular point}. We denote the set of regular points by $\mathcal{R}(u)$ and the set of singular points by $\mathcal{S}(u)$.
\end{definition}

If $p$ is a regular point of $u$, then locally $u$ is equal to the composition of a harmonic map to a closed interval $I\subset \mathbb{R}$, and an isometric embedding of $I$ in $\mathcal{T}$. It is clear from the definition that $\mathcal{R}(u)$ is open and $\mathcal{S}(u)$ is closed.

\begin{theorem}[{\cite[Thm.\ 1.4]{sun2003regularity}}]
	\label{thm_sun_locally_finite}
	Let $u: \Omega \to \mathcal{T}$ be a harmonic map from a manifold to an $\mathbb{R}$--tree. Then, for every point $p \in \Omega$, there exists an open ball $B_r(p)$ centered at $p$ such that  $u(B_{r}(p))$ lies in an embedded locally finite subtree.
\end{theorem}

Theorem \ref{thm_sun_locally_finite} allows one to reduce the regularity
problem of harmonic maps into $\mathbb{R}$--trees to the case where the
tree is locally finite. Before introducing the next results, we review the
definition of the order function from \cite[Sec.\ 2]{gromov1992harmonic}.

\begin{definition}
	Assume $u:\Omega\to \mathcal{T}$ is harmonic. 
	For $p \in \Omega$, the \emph{order} of $u$ at $p$ is defined to be the limit 
	\begin{equation}
		\label{eqn_def_ord}
			\mathrm{ord}(u; p) := \lim_{r \to 0}  \frac{r \int_{B_r(p)} |\nabla u(x)|^2\,\mu(x)}{\int_{\partial B_r(p)} d^2(u(p), u(x))\,\sigma(x)},
	\end{equation}
	where $\mu(x)$ denotes the volume measure of $\Omega$ and $\sigma(x)$ denotes the area form.
\end{definition}
It was proved in \cite[Sec.\ 2]{gromov1992harmonic} that the limit in \eqref{eqn_def_ord} always exists. 
When $\mathcal{T}=\mathbb{R}$, the value of $\mathrm{ord}(u; p)$ equals the
vanishing order of $u-u(p)$ at $p$. Moreover, the order function
$\mathrm{ord}(u;p)$ is upper semi-continuous with respect to $p$
(\cite[statement above Thm.\ 2.3]{gromov1992harmonic}).

We collect the following regularity results from the literature, which will be used later. In the following,  $u:\Omega\to \mathcal{T}$ denotes a harmonic map, and $d$ denotes the metric on the $\mathbb{R}$--tree $\mathcal{T}$. 

\begin{theorem}[{\cite[Thm.\ 6.3]{gromov1992harmonic}}, {\cite[Thm.\ 1.1]{sun2003regularity}}]
	\label{thm_epsilon_regularity}
	There is a constant $\epsilon$ depending only on the dimension of the domain, such that for each point $p \in \Omega$, either $\mathrm{ord}(u; p) =1$, or $\mathrm{ord}(u; p) \ge 1+\epsilon$. Moreover, if $\mathrm{ord}(u; p) = 1$, then $p$ is a regular point. 
\end{theorem}

\begin{theorem}[{\cite[Thm.\ 2.4.6]{korevaar1993sobolev}}] 
	\label{thm_harmonic_lipschitz}
	The map $u$ is locally Lipschitz, and the Lipschitz constants only depend on the curvature of $\Omega$. 
\end{theorem}

\begin{theorem}
	[{\cite[Thm.\ 1.3]{sun2003regularity}}, {\cite[Thm.\ 6.4]{gromov1992harmonic}}]
	\label{thm_harmonic_hausdorff_dim}
	The singular set $\mathcal{S}(u)$ has Hausdorff codimension at least $2$. 
\end{theorem}

A stronger version of this result was given by Dees \cite{dees2022rectifiability}.

\begin{theorem}
	[{\cite[Thms.\ 1.1 and 1.2]{dees2022rectifiability}}]
	\label{thm_harmonic_rectifiable}
	Assume $\Omega$ has dimension $n$, then $\mathcal{S}(u)$ is $(n-2)$--rectifiable. Moreover, for each compact set $K\subset \Omega$, there exists a constant $C$ depending on $K$, $u$, and $\Omega$, such that 
	$$
	\mathrm{Vol} (B_r(\mathcal{S}(u)\cap K)) \le C \cdot r^2
	$$
	for $r<1$, 
	where $B_r(\mathcal{S}(u)\cap K)$ denotes $r$--neighborhood of $\mathcal{S}(u)\cap K$. 
\end{theorem}

\begin{theorem}[{\cite[proof of Thm.\ 2.3]{gromov1992harmonic}}] 
	\label{thm_holder_estimate_from_frequency}
	Assume $\Omega$ as a manifold is the open unit ball $B_1(0)$ in $\mathbb{R}^n$, but the metric is not necessarily Euclidean. Suppose $\mathrm{ord}(u; 0) = \alpha$. Then there exists a constant $C$ depending only on the curvature, such that for all $\lambda\in[0,1]$ and $x\in B_{1/2}(0)$, we have
	\begin{equation}
		\label{eq_homogeneous_estimate}
		d(u(\lambda x), u(0)) \leq C \lambda^{\alpha} d(u(x), u(0)).
	\end{equation}
\end{theorem}

\begin{corollary}
	\label{cor_L2_Holder_bound}
	Assume $\mathrm{ord}(u;p)=\alpha$. Let $n$ be the dimension of $\Omega$. Let $r_0$ be the minimum of $1$ and the injectivity radius of $\Omega$ at $p$.  Assume $r<r_0/2$. Then there exists a constant $C$ depending only on the curvature of $\Omega$ in $B_{r_0}(p)$, such that 
	$$
	\int_{B_r(p)} |\nabla u|^2 \le C r^{n-2+2\alpha}\cdot r_0^{2-2\alpha}.
	$$
\end{corollary}

\begin{proof}
	By  \cite[Sec.\ 2]{gromov1992harmonic}, there exists a constant $C_1$ depending only on the curvature, such that 
	$$
	e^{C_1r^2}\frac{r \int_{B_r(p)} |\nabla u|^2}{\int_{\partial B_r(p)} d^2(u(p), u)}
	$$
	is increasing with respect to $r$ for all $r<r_0$. 
	Hence 
	$$
	e^{C_1r^2}\frac{r \int_{B_r(p)} |\nabla u|^2}{\int_{\partial B_r(p)} d^2(u(p), u)} \le e^{C_1r_0^2}\frac{r_0 \int_{B_{r_0}(p)} |\nabla u|^2}{\int_{\partial B_{r_0}(p)} d^2(u(p), u)}.
	$$
	By Theorem \ref{thm_harmonic_lipschitz}, we have $r_0 \int_{B_{r_0}(p)} |\nabla u|^2\le C_2{r_0}^{n+1}$ for a constant $C_2$. 
	By Theorem \ref{thm_holder_estimate_from_frequency}, we have 
	$$
	\frac{\int_{\partial B_r(p)} d^2(u(p), u)}{\int_{\partial B_{r_0}(p)} d^2(u(p), u)}\le \left(\frac{r}{r_0}\right)^{n-1} \left(\frac{r}{r_0}\right)^{2\alpha},
	$$
	where the first factor comes from the comparison of area forms. The desired result then follows from a straightforward computation. 
\end{proof}

We also need the following estimate.

\begin{theorem}[{\cite[Thm.\ 2.4]{gromov1992harmonic}}]
	\label{thm_C0_from_L2_harmonic}
	Assume $p\in \Omega$ and $r$ is less than the injectivity radius of $\Omega$ at $p$. Then there exists a constant $C$, depending only on the curvature of $\Omega$ on $B_r(p)$, such that 
	$$
	\|\nabla u\|_{L^\infty(B_{r/2}(p))}^2 \le C \fint_{B_r(p)} |\nabla u|^2.
	$$
\end{theorem}

\begin{remark}
	\label{rmk_nabla_u_continuous}
	Note that the definitions in Section \ref{subsec_energy_density} only defined $|\nabla u|$ almost everywhere on $\Omega$. On $\mathcal{R}(u)$, the function $|\nabla u|$ is represented by a smooth function. By Theorem \ref{thm_C0_from_L2_harmonic} and Corollary \ref{cor_L2_Holder_bound}, we see that the function $|\nabla u|$ is given by a continuous function on $\Omega$, where $|\nabla u| = 0$ on $\mathcal{S}(u)$.
\end{remark}

We now prove the following estimate.
\begin{proposition}
	\label{prop_nabla_nabla_u_finite}
	Assume $K$ is a compact subset of $\Omega$. Then 
	$$
	\int_{K\setminus \mathcal{S}(u)} |\nabla \nabla u|^2 <+\infty.
	$$
\end{proposition}

\begin{proof}
	By Theorem \ref{thm_harmonic_hausdorff_dim}, $\mathcal{S}(u)$ has zero Lebesgue measure. Hence we may  write $\int_{\Omega\setminus \mathcal{S}(u)}$ as $\int_\Omega$, if the integrand is defined on $\mathcal{R}(u)$ and extends to an integrable function on $\Omega$.
	
	By shrinking $\Omega$ to a smaller domain containing $K$, we may assume that the curvature of $\Omega$ is uniformly bounded and $\Omega$ has finite volume. By Theorem \ref{thm_harmonic_lipschitz}, this implies $|\nabla u|$ is bounded and $\int_\Omega |\nabla u|^2 <+\infty$. 
	
	By \cite[eq.\ (6.2)]{gromov1992harmonic}, there is a constant $C$ depending only on the curvature of $\Omega$ such that 
	\begin{equation}
		\label{eqn_harmonic_func_laplacian}
		\frac12 \Delta|\nabla u|^2 \ge |\nabla\nabla u|^2 - C|\nabla u|^2 
	\end{equation}
	pointwise on $\mathcal{R}(u)$.

	Let $\rho$ be a Lipschitz function that is compactly supported in $\mathcal{R}(u)$ such that $0\le \rho \le 1$. By \eqref{eqn_harmonic_func_laplacian}, 
	$$
	\frac12 \int_\Omega \Delta|\nabla u|^2\rho^2 \ge \int_\Omega |\nabla \nabla u|^2 \rho^2 - C\int_\Omega |\nabla u|^2 \rho^2,
	$$
	so
	\begin{align*}
		\int_\Omega |\nabla \nabla u|^2 \rho^2
		& \le C\int_\Omega |\nabla u|^2 \rho^2 + \frac12 \int_\Omega \Delta|\nabla u|^2\rho^2 \\
		& = C\int_\Omega |\nabla u|^2 \rho^2 - \frac12 \int_\Omega \langle \nabla|\nabla u|^2 ,\nabla( \rho^2) \rangle \\
		& \le C\int_\Omega |\nabla u|^2 \rho^2 + 2\int_\Omega |\nabla|\nabla u||\cdot |\nabla u| \cdot |\nabla \rho|\cdot \rho \\
		& \le C\int_\Omega |\nabla u|^2 \rho^2 + 2\int_\Omega |\nabla\nabla u|\cdot |\nabla u| \cdot |\nabla \rho|\cdot \rho \\
		& \le C\int_\Omega |\nabla u|^2 \rho^2 + \frac12 \int_\Omega |\nabla \nabla u|^2 \rho^2 + 2\int_\Omega |\nabla u|^2 |\nabla \rho|^2 .\\
	\end{align*}
	As a result, 
	\begin{equation}
		\label{eqn_nabla_nabla_u_estimate}
		\int_\Omega |\nabla \nabla u|^2 \rho^2 \le 2 C\int_\Omega |\nabla u|^2 \rho^2 + 4 \int_\Omega |\nabla u|^2 |\nabla \rho|^2.
	\end{equation}
	By the assumptions on $\Omega$, we have 
	$$
	2 C\int_\Omega |\nabla u|^2 \rho^2 \le 2C\|\nabla u\|_{L^\infty}^2 \,\mathrm{Vol}(\Omega).
	$$
	
	For every positive integer $i$, let $\nu_i:[0,+\infty)\to [0,1]$ be the Lipschitz function such that $\nu(x) = 1$ for $x\ge 2/i$, $\nu(x)=0$ for $x\le 1/i$, and $\nu'(x) = i$ for $x\in [1/i,2/i]$. Let $\chi$ be a smooth function that is compactly support in $\Omega$, such that $0\le \chi\le 1$, and $\chi=1$ on the given compact set $K$. Define $\rho_i$ by 
	$$
	\rho_i(x) = [\nu_i\circ d_\Omega (K\cap \mathcal{S}(u),x)]\cdot \chi,
	$$
	where $d_\Omega$ is the distance function on $\Omega$.  Theorem \ref{thm_harmonic_rectifiable} implies that $\int_\Omega |\nabla \rho_i|^2$ is bounded as $i\to\infty$. Hence by \eqref{eqn_nabla_nabla_u_estimate}, the integral $\int_\Omega |\nabla \nabla u|^2\rho_i^2$ is bounded as $i\to\infty$. Since $\rho_i$ is increasing with respect to $i$ and converges to $\chi$ as $i\to \infty$, this implies $\int_\Omega \chi\cdot |\nabla u|^2$ is finite, so the desired result is proved. 
\end{proof}

\subsection{From harmonic maps to $\ZT$ harmonic 1-forms}
\label{subsec_harmonic_map_to_form}
Now assume $M$ is a closed oriented $3$--manifold and let $\widetilde{M}$ be its universal cover. Let 
$
\pi:\widetilde{M}\to M
$
be the covering map.
Assume $(\mathcal{T},d_{\mathcal{T}})$ is an $\mathbb{R}$--tree with an isometric $\pi_1(M)$--action, and assume $u:\widetilde{M} \to \mathcal{T}$ is a $\pi_1(M)$--equivariant harmonic map. 

Then $\nabla u$ defines a two-valued $1$--form on $\mathcal{R}(u)\subset \widetilde{M}$. It is two-valued because the geodesics on $\mathcal{T}$ do not have canonical orientations. Locally, $\nabla u$ is given by $\pm v$ for some (single-valued) harmonic $1$--form $v$. Since $u$ is $\pi_1(M)$ equivariant, $\nabla u$ defines a two-valued $1$--form on the image of $\mathcal{R}(u)$ in $M$.

\begin{theorem}
	\label{thm_construction_Z2harmonic}
	There exists a unique $\ZT$ harmonic $1$--form $\bv$ on $M$, such that $\pi^*(\bv) = \nabla u$ on $\mathcal{R}(u)$. 
\end{theorem}

\begin{proof}
	The uniqueness is clear since $\mathcal{S}(u)$ has Hausdorff codimension $2$.  We prove the existence of such $\bv$.
	Recall that $ |\nabla u|$ is continuous (see Remark \ref{rmk_nabla_u_continuous}). 
	Let $\widetilde{Z} = |\nabla u|^{-1}(0)$. Then $\mathcal{S}(u)\subset \widetilde{Z}$, so $\nabla u$ defines a non-vanishing two-valued $1$--form $\tilde\bv$ on $\widetilde{M}\setminus \widetilde{Z}$. Since $u$ is equivariant, the set $\widetilde{Z}$ is invariant under the $\pi_1(M)$--action. Let $Z$ be the image of $\widetilde{Z}$ in $M$. Then $\tilde\bv$ is the pull-back of a two valued $1$--form $\bv$ on $M\setminus Z$. 
	
	We show that $\bv$ is a $\ZT$ harmonic $1$--form with respect to the zero locus $Z$. Since $\nabla u$ is harmonic in $\mathcal{R}(u)$, we know that $\bv$ is locally given by non-vanishing harmonic $1$--forms on $M\setminus Z$.
	
	Note that if $p\in \mathcal{R}(u)$ and $|\nabla u (p)|= 0$, then $\mathrm{ord}(u;p)\ge 2$. As a result, by Theorem \ref{thm_epsilon_regularity}, there exists $\epsilon>0$ such that $\mathrm{ord}(u;p)\ge 1+\epsilon$ for all $p\in \widetilde{Z}$. By Corollary \ref{cor_L2_Holder_bound}, we have
	$$
	\int_{B_r(p)} |\nabla u|^2 \le C r^{3+2\epsilon}\cdot r_0^{2-2\mathrm{ord}(u;p)}
	$$
	for all $p\in \widetilde{Z}$ and $r<r_0/2$, where $r_0$ is the minimum of $1$ and the injectivity radius of $M$ at $p$. Since $\mathrm{ord}(u;p)$ is upper semi-continuous and $M$ is closed, it is bounded from above on $\widetilde{M}$. Since $M$ is closed, the value of $r_0$ has a positive lower bound. Therefore, the value of $r_0^{2-2\mathrm{ord}(u;p)}$ has an upper bound (which may depend on $u$). This verifies Condition (iii) of Definition \ref{def_ZTharmonicform}.
	
	Condition (ii) of Definition \ref{def_ZTharmonicform} follows immediately from the fact that $|\nabla u|$ is  bounded and Proposition \ref{prop_nabla_nabla_u_finite}.
\end{proof}

\begin{remark}
	\label{rmk_relation_Z_R(u)}
	From the above construction, we also see that $\mathcal{S}(u)\subset \widetilde{Z}$. In other words, every non-vanishing point of the $\ZT$ harmonic $1$--form corresponds to regular points of $u$ on $\widetilde{M}$. 
\end{remark}

\begin{definition}
\label{def_associated_ZT_spinor}
	We call the $\ZT$ harmonic form $\bv$ obtained by Theorem \ref{thm_construction_Z2harmonic}  the $\ZT$ harmonic 1-form \emph{associated with $u$}.
\end{definition}

\subsection{Maps between trees}
Let $M$, $\widetilde{M}$, $(\mathcal{T},d)$, $u$, be as in Section \ref{subsec_harmonic_map_to_form}. Let $\bv$ be the $\ZT$ harmonic $1$--form on $M$ associated with $u$. Let $\tilde{\bv}$ be the pull back of $\bv$ to $\widetilde{M}$. By Theorem \ref{thm_pseudo-metric-R-tree}, the leaf space $(\mathcal{T}_{\widetilde{M},\tilde{\bv}}, d_{\widetilde{M},\tilde{\bv}})$  is also an $\mathbb{R}$--tree. In this subsection, we study the relationship between $\mathcal{T}_{\widetilde{M},\tilde{\bv}}$ and $\mathcal{T}$. 

Let $\mu_{\tilde \bv}$ be the transverse measure defined from $\tilde \bv$ by \eqref{eqn_transverse_measure_from_bv}. Let $Z\subset M$, $\widetilde{Z}\subset \widetilde{M}$ denote the zero loci of $\bv$ and $\tilde\bv$. 

\begin{lemma}
	\label{lem_compare_measure_and_image}
	Suppose $\gamma:[0,1]\to \widetilde{M}$ is a $\MC^1$ arc. Let $p=\gamma(0)$, $q=\gamma(1)$. Then 
	\begin{equation}
		\label{eqn_compare_mu_d}
		\mu_{\tilde\bv}(\gamma) \ge d(u(p),u(q)).
	\end{equation}
	Moreover, equality holds if $\gamma$ is transverse to $\ker \tilde\bv$. 
\end{lemma}
\begin{proof}
	We first prove \eqref{eqn_compare_mu_d}.
	Both sides of \eqref{eqn_compare_mu_d} are continuous with respect to $\gamma$ in the $\MC^1$ topology, so the inequality is a closed condition. Recall that  $\mathcal{R}(u)\subset \widetilde{M}$ denotes the regular set of $u$.  Since $\mathcal{R}(u)$ has Hausdorff codimension at least $2$, after perturbing $\gamma$ if necessary, we may assume without loss of generality that the image of $\gamma$ is contained in $\mathcal{R}(u)$. 
	
	Then we have 
	\begin{equation}
		\label{eqn_compare_mu_d_derivative}
		\frac{d}{dt} d\bigg(u(p),u\big(\gamma(t)\big)\bigg) 
		\le\big |\big \lan\nabla u(\gamma(t)), \dot\gamma(t)\big \ran\big|\\
		=\frac{d}{dt} \mu_{\tilde\bv} (\gamma|_{[0,t]}).
	\end{equation}
	Hence \eqref{eqn_compare_mu_d} is proved by integrating \eqref{eqn_compare_mu_d_derivative} for $t$ in $[0,1]$. 
	
	If $\gamma$ is transverse to $\ker \tilde\bv$, then the image of $\gamma$ is contained in the complement of $\widetilde Z$. By Remark \ref{rmk_relation_Z_R(u)}, the image of $\gamma$ is contained in $\mathcal{R}(u)$. Then for $t\in [0,1]$, the value of $u\circ \gamma(t)\in\MT$ locally moves along geodesics at non-vanishing velocities. Since $\MT$ is an $\mathbb{R}$--tree, this implies that the image of $u\circ \gamma(t)$ for $t\in[0,1]$ is a geodesic segment where the point $u\circ \gamma(t)$ moves at non-vanishing velocities with respect to $t$. Hence the inequality in \eqref{eqn_compare_mu_d_derivative} achieves equality for all $t\in [0,1]$, and the desired result is proved. 
\end{proof}

Let $\pi_{\tilde{\bv}}$ denote the quotient map from $\widetilde{M}$ to $\MT_{\widetilde{M}, \tilde{\bv}}$.

\begin{theorem}
	\label{thm_factorize_harm_map_thru_leafs}
	There exists a unique continuous map $f: \mathcal{T}_{\widetilde{M}, \tilde{\bv}} \to \mathcal{T}$ such that $u = f \circ \pi_{\tilde{\bv}}$. Moreover, we have the following properties:
	\begin{enumerate}
		\item $f$ is $1$--Lipschitz.
		\item The map $\pi_{\tilde{\bv}}$ is harmonic.
		\item $|\nabla \pi_{\tilde{\bv}}| = |\nabla u|$.
		\item If $|\nabla u(p)| \neq 0$, then $\nabla \pi_{\tilde{\bv}}(p) = \nabla u(p)$ as two-valued harmonic 1-forms.  
	\end{enumerate}
\end{theorem}

\begin{proof}
	By Lemma \ref{lem_compare_measure_and_image}, we have 
	$$
	d_{\widetilde{M},\tilde{\bv}}(\pi_{\tilde\bv}(p), \pi_{\tilde\bv}(q)) \ge d(u(p),u(q))
	$$
	for all $p,q\in \widetilde{M}$. Hence the map $u$ factorizes uniquely as $u=  f\circ \pi_{\tilde\bv}$ for a continuous map $f$, and $f$ is $1$--Lipschitz. As a result, we have  $|\nabla u|\le |\nabla \pi_{\tilde\bv}|$.
	
	By the definition of $d_{\widetilde{M},\tilde{\bv}}$, the map $\pi_{\tilde\bv}$ is Lipschitz with Lipschitz constant $\|\bv\|_{\MC^0}$, so the map $\pi_{\tilde\bv}$ is in $W^{1,2}_{loc}$.
	
	For $p\in \widetilde{M}\setminus \widetilde{Z}$, let $\gamma:[0,1]\to \widetilde{M}$ be an arc that is transverse to $\ker \tilde\bv$.
	Then 
	$$
	d_{\widetilde{M},\tilde{\bv}}\big(\pi_{\tilde\bv}\circ \gamma(0), \pi_{\tilde\bv}\circ \gamma(1)\big) \ge d(u(\gamma(0)),u(\gamma(1))) = \mu_{\tilde\bv}(\gamma) \ge 	d_{\widetilde{M},\tilde{\bv}}\big(\pi_{\tilde\bv}\circ \gamma(0), \pi_{\tilde\bv}\circ \gamma(1)\big).
	$$
	Here, the first inequality follows from the fact that $f$ is $1$--Lipschitz, the second equality follows from Lemma \ref{lem_compare_measure_and_image}, and the third inequality follows from the definition of $d_{\widetilde{M},\tilde{\bv}}$. Hence $	d_{\widetilde{M},\tilde{\bv}}\big(\pi_{\tilde\bv}\circ \gamma(0), \pi_{\tilde\bv}\circ \gamma(1)\big) = d(u(\gamma(0)),u(\gamma(1)))$. On the other hand, if $\gamma$ is tangent to $\ker \tilde\bv$ and is contained in the complement of $\widetilde{Z}$, then both $	d_{\widetilde{M},\tilde{\bv}}\big(\pi_{\tilde\bv}\circ \gamma(0), \pi_{\tilde\bv}\circ \gamma(1)\big)$ and $d(u(\gamma(0),u(\gamma(1)))$ are zero. In conclusion, if $p$ is in the complement of $\widetilde{Z}$, then there exists an open ball $B_r(p)$ centered at $p$ such that 
	$$
	d_{\widetilde{M},\tilde{\bv}}\big(\pi_{\tilde\bv} (p), \pi_{\tilde\bv} (q)\big)
	= d(u(p),u(q))
	$$
	for all $q\in B_r(p)$.
	This implies that the directional derivatives and energy densities of $\pi_{\tilde\bv}$ and $u$ are the same on $\widetilde{M}\setminus \widetilde{Z}$. 
	Since $\widetilde{Z}$ has Lebesgue measure zero, they also define the same measure density functions on $\widetilde{M}$. 
	
	It remains to show that $\pi_{\tilde\bv}$ is harmonic. Let $M_0\subset \tM$ be a fundamental domain of the $\pi_1$--action. An equivariant map from $\widetilde{M}$ to a tree is harmonic if and only if it minimizes the energy on $M_0$ among all locally $W^{1,2}$ equivariant maps. We show that $\pi_{\tilde\bv}$ minimizes the energy on $M_0$.  Assume there exists equivariant map $g:\widetilde M \to \mathcal{T}_{\widetilde{M},\tilde\bv}$ such that 
	$$
	\int_{M_0}|\nabla g|^2 < \int_{M_0}|\nabla \pi_{\tilde\bv}|^2,
	$$
	then we have
	$$
	\int_{M_0}|\nabla (f\circ g)|^2 \le \int_{M_0}|\nabla  g|^2 < \int_{M_0}|\nabla \pi_{\tilde\bv}|^2 =\int_{M_0}|\nabla u|^2,
	$$
	where the first inequality follows from the fact that $f$ is $1$--Lipschitz, and the last equation follows from the fact that  $|\nabla \pi_{\tilde\bv}| =|\nabla u|$ on $\widetilde{M}$. 
	Since $f$ is $\pi_1(M)$--equivariant, this contradicts the assumption that $u$ is energy minimizing. Hence the theorem is proved. 
\end{proof}

\begin{remark}
	In general, if $(M,\bv)$ satisfies the conditions of Theorem \ref{thm_pseudo-metric-R-tree}, the quotient map from $M$ to the leaf tree $\MT_{M,\bv}$ may not necessarily be harmonic. See Section \ref{subsec_proj_leaf_space} for a counterexample.
\end{remark}

\section{Proofs of the main theorems}
In this section, we set up the relevant constructions and prove Theorem
\ref{thm_map_on_closure}. We will also show that Theorem
\ref{thm_factoriztion_tree_maps} follows directly from Theorem
\ref{thm_factorize_harm_map_thru_leafs}, once the map $\MH$ in \eqref{eq_boundarymap} is constructed. Following our previous convention, let $M$ denote an oriented closed Riemannian 3-manifold, and let $\widetilde{M}$ be its universal cover. Let $\Gamma = \pi_1(M)$.

\subsection{Definition of the boundary map $\MH$}
\label{subsection_construction_boundary_map}
 Recall that $\MPL(\Gamma)$ denotes the space of
projective length functions of all complete minimal $\Gamma$--trees, and
$\MMZ$ denotes the space of all $\ZT$ harmonic $1$--forms on $M$ with unit $L^2$--norm. 
In this subsection, we define a map 
\begin{equation}
	\label{eq_boundarymapH}
\MH: \MPL(\Gamma) \to \MMZ.	
\end{equation}

By Theorem \ref{thm_nonabelian_determinant_tree}, if $[\ell]\in
\MPL(\Gamma)$ is non-abelian, then there is a unique (up to equivariant isometry) minimal
$\Gamma$--tree $\mathcal{T}$ with length function  $\ell$ and no fixed ends. 
By Theorem \ref{thm_existence_uniqueness} there is a $\pi_1(M)$ equivariant harmonic map $u:\widetilde{M}\to\MT$. By Proposition \ref{prop_uniqueness_directional_energy}, (see also Definition \ref{def_associated_ZT_spinor}), the $\ZT$ harmonic $1$--form associated with $u$ is independent of the choice of $\ell$ or $u$. Let $\bv$ be the $\ZT$ harmonic $1$--form associated with $u$. Define 
$$
\MH([\ell]) = \bv/\|\bv\|_{L^2(M)}.
$$
Since $\bv/\|\bv\|_{L^2(M)}$ only depends on the projective class of $\ell$ in $\MPL(\Gamma)$, the value of $\MH([\ell])$ is well-defined. 

Next, we define the map $\MH$ on abelian length functions. Assume $\ell$ is an abelian length function.
By definition, there exists a homomorphism $\mu:\Gamma \to \mathbb{R}$ so that $\ell(x) = |\mu(x)|$ for all $x\in \Gamma$. It follows from straightforward algebra that if $\mu_1,\mu_2$ are two homomorphisms
$\Gamma\to \mathbb{R}$ such that $|\mu_1(x)| = |\mu_2(x)|$ for all $x\in \Gamma$, then $\mu_1=\pm \mu_2$. In other words, $\mu$ is determined by $\ell$ up to an overall sign. 

Define $v_\mu$ to be the (single-valued) harmonic $1$--form on $M$ such that the periods of $v_\mu$ on $\pi_1(M)$ are given by $\mu$. Let $\bv$ be the $\ZT$ harmonic $1$--form given by $\pm v_\mu$. Define  
$$
\MH([\ell]) = \bv/\|\bv\|_{L^2(M)}.
$$
In conclusion, we have constructed a well-defined map $\MH$ from $\MPL(\Gamma)$ to  $\MMZ$.

\begin{lemma}
	\label{lem_abelian_length}
Suppose $\MT$ is a minimal $\Gamma$ tree whose length function $\ell$ is abelian and not identically zero, and $u:\widetilde{M}\to \MT$ is a $\Gamma$--equivariant harmonic map. Then the tree  $\mathcal{T}$ must be $\mathbb{R}$, and the $\Gamma$ action on $\mathcal{T}$ is given by  translations.
\end{lemma}

\begin{proof}
	Suppose that $\Gamma$ has a minimal action on a tree $\mathcal{T}$ not
    isometric to $\mathbb{R}$, but with an abelian length function. By
    \cite[p.\ 573]{cullermorgan97group}, the action of $\Gamma$ must have a
    fixed  end. Combining \cite[Thm.\ 5.3]{daskalopoulos1998character}
    and \cite[Thm.\ 1.2]{mese2002uniqueness}, we see that there can be no $\Gamma$-equivariant harmonic map to $\mathcal{T}$. 
	
	Now that we have proved $\MT\cong\mathbb{R}$, we show that the $\Gamma$ action must be given by translations. Assume $x\in\Gamma$ acts on $\MT\cong \mathbb{R}$ by an orientation-reversing isometry, then $\ell(x)=0$. As a result, for every $y\in \Gamma$, we must have $\ell(xy) = \ell(y)$. This is only possible when the action of $y$ equals the identity or the action of $x$. Hence the length function of the $\Gamma$--action is identically zero, contradicting the assumptions.
\end{proof}

\begin{corollary}
	\label{cor_MH(ell)_distance_function}
	Assume $\mathcal{T}$ is a complete minimal $\Gamma$--tree and $u:\widetilde{M}\to \MT$ is a $\Gamma$--equivariant harmonic map. Let $\ell$ be the length function of $\MT$, and assume $\ell$ is not identically zero. Then the $\ZT$ harmonic $1$--form associated with $u$ is equal to $\MH([\ell])$ up to a non-zero constant multiplication.
\end{corollary}

\begin{proof}
	If $\ell$ is non-abelian, the statement follows from the definition of $\MH$.
	If $\ell$ is abelian, then Lemma \ref{lem_abelian_length} implies that $\MT\cong \mathbb{R}$ and the action of $\Gamma$ on $\MT$ are translations. Therefore, $u$ is a harmonic function and the associated $\ZT$ harmonic $1$--form is given by $\pm du$. Since the periods of $du$ coincide with the translation distances of the $\Gamma$--action on $\MT$, the desired statement is proved. 
\end{proof}

\subsection{Comparison of the compactification limits}
Recall that $\cX(\Gamma)$ denotes the character variety of $\Gamma$, and that $\MMc$ denotes the moduli space of solutions to \eqref{eqn_KW_system}. Also recall that $\cX(\Gamma)$  and $\MMc$  are canonically identified by the Riemann--Hilbert correspondence. In this subsection, we prove the following proposition.

\begin{proposition}
	\label{prop_compare_two_limits}
	Assume $\{p_i\}$ is a sequence of points in $\cX(\Gamma)$, and let $\{p_i'\}$ be the corresponding sequence in $\MMc$. Assume that $p_i$ converges in the Morgan--Shalen compactification to the projective length function $[\ell]$, and that $p_i'$ converges in the Taubes compactification $\BMMc$ to a $\ZT$ harmonic $1$--form $\bv$. Then $\MH([\ell]) = \bv$. 
\end{proposition}

The proof of Proposition \ref{prop_compare_two_limits} amounts to careful book-keeping using results from the earlier sections. We start by reviewing the characterizations of the Morgan--Shalen compactification and the Taubes compactification.

\subsubsection {Characterization of the Morgan--Shalen limit}
\label{subsubsec_MS_lim}
By \cite{daskalopoulos1998character}, the Morgan--Shalen limit is given by the limit of harmonic maps from $\widetilde{M}$ to $\mathbb{H}^3$. More precisely, $\mathbb{H}^3$ denotes the three-dimensional hyperbolic space form. Its orientation-preserving isometry group is isomorphic to $\mathrm{PSL}_2(\mathbb{C})$. Every $p_i\in \chi(\Gamma)$ is given by a representation of $\Gamma$ in $\SLC$, which defines a $\Gamma$--action on $\mathbb{H}^3$. Let $u_i$ be a $\Gamma$--equivariant harmonic map from $\widetilde{M}$ to $\mathbb{H}^3$ with respect to the above $\Gamma$--action. Let $E(u_i)$ be the energy of $u_i$ on a fundamental domain, and let $d_i$ denote the pseudo-metric on $\widetilde{M}$ defined by 
\begin{equation}
	\label{eqn_d_i}
d_i(p,q) = \frac{d_{\mathbb{H}^3}(u_i(p),u_i(q))}{ E(u_i)^{1/2}}.
\end{equation}
If $\{p_i\}$ converges to a projective length function $[\ell]$ on the boundary of the Morgan--Shalen compactification, then $E(u_i)\to \infty$, and $d_i$ converges locally uniformly to a limit $d_\infty$.  The quotient metric space of $\widetilde{M}$ with respect to $d_\infty$ is a minimal $\Gamma$--tree $\MT$, the length function of $\MT$ is in the projective class $[\ell]$, and the quotient map from $\widetilde{M}$ to $\MT$ is harmonic. 

\subsubsection{Characterization of the Taubes limit}
If $\{p_i'\}$ converges to a $\ZT$ harmonic $1$--form $\bv$, assume each $p_i'$ is given by the pair $(A_i,\phi_i)$. Let $P_i$ be the principal $\SU(2)$ bundle over $M$ where $A_i$ is defined. Let $\mathfrak{g}_{P_i}$ be the associated $\mathfrak{su}(2)$ bundle with respect to the adjoint action. 

Let $Z$ be the zero locus of $\bv$. By the constructions in Section \ref{sec_Z2_forms}, on every compact set $K$ that is contained in an open ball in $M\setminus Z$, there exists $\phi$ on $K$, such that after a sequence of gauge transformations, $\phi_i/\|\phi_i\|_{L^2(M)}$ converge to $\phi$ in the weak $W^{2,2}$ topology on $K$. By Sobolev embedding theorems,  $\phi_i/\|\phi_i\|_{L^2(M)}$ converge to $\phi$ in the $C^0$ topology. 

By part (3) of the second bullet point of \cite[Thm.\ 1.1a]{Taubes20133manifoldcompactness}, the limit spinor $\phi$ locally has the form $\sigma\otimes v$, where $v$ is a harmonic $1$--form, $\sigma$ is a section of $\mathfrak{g}_{P_i}$, and $|\sigma|=1$ pointwise. In this case, the $\ZT$ harmonic $1$--form $\bv$ is locally equal to $\pm v$. 

\subsubsection{Relationship between the two constructions}
\label{subsubsec_relation_two_lim}
The relationship of the maps $u_i$ and $\phi_i$ is given as follows (see
\cite[Sec.\ 2C]{daskalopoulos1998character}). 
Let $\widetilde{P}_i$ and $\tilde{\mathfrak{g}}_{P_i}$ denote the pull-backs of $P_i$ and $\mathfrak{g}_{P_i}$ to $\widetilde{M}$. 
The orthonormal frame bundle of $\mathbb{H}^3$ is an $\SO(3)$ bundle. Its Lie algebra bundle is an $\mathfrak{so}(3)$ bundle, and the pull-back of this bundle to $\widetilde{M}$ is isomorphic to  $\tilde{\mathfrak{g}}_{P_i}$. Then there exists an equivariant isomorphism such that $-\frac12 \nabla u_i =  \phi_i$ under the identification of Lie algebra bundles, where $\nabla$ is given by the Levi-Civita connections of $\widetilde{M}$ and $\mathbb{H}$. Different choices of equivariant isomorphisms of Lie algebra bundles correspond to gauge transformations on solutions to \eqref{eqn_KW_system}.

As a consequence, we have 
\begin{equation}
	\label{eqn_relation_u_phi}
\frac12 |\nabla u_i(X)| = |\phi_i(X)|
\end{equation}
for every tangent vector $X$ of $\widetilde{M}$. Note that this equation is independent of the choice of the isomorphism between Lie algebra bundles, or gauge transformations. 

\subsubsection{Proof of Proposition \ref{prop_compare_two_limits}}
Let $u:\widetilde{M}\to \MT$ denote the quotient map from $\widetilde{M}$ to $\MT$. Recall that $\bv$ denotes the $\mathbb{Z}/2$ harmonic form in the Taubes limit. Let $\tilde \bv$ denote the pull-back of $\bv$ to $\widetilde{M}$.
Let $\mathcal{U}$ be the set of $x\in \widetilde{M}$ such that $x$ is not on the zero locus of $\tilde \bv$ or the zero set of $|\nabla u|$.
Then $\mathcal{U}$ is a  $\Gamma$--equivariant open and dense subset of $\widetilde{M}$. The set $\mathcal{U}$ is contained in the regular set of $u$, so $\nabla u$ is a non-vanishing $\ZT$ harmonic $1$--form on $\mathcal{U}$. 
 
As a result, $\ker \bv$ and $\ker \nabla u$ are two smooth foliations on $\mathcal{U}$. 
We first show that these two foliation are the same. 

For each $p\in \mathcal{U}$, let $B_r(p)$ be an open ball centered at $p$ such that 
\begin{enumerate}
	\item $B_r(p)\subset \mathcal{U}$, 
	\item the metrics $d_i$ (defined in Equation \eqref{eqn_d_i}) converge to $d_\infty$ uniformly on $B_r(p)$, 
	\item $\phi_i$ uniformly converge to a limit $\phi$  on $B_r(p)$ after gauge transformations.
\end{enumerate}

 Let $\gamma:[0,1]\to B_r(p)$ be a smooth arc that is tangent to $\ker \bv$. Write $x = \gamma(0)$, $y=\gamma(1)$. By \eqref{eqn_relation_u_phi}, we have
\begin{align*}
d_\infty(x,y) & = \lim_{i\to\infty} d_i(x,y) \le \limsup_{i\to \infty} \int_0^1 |\nabla u_i(\dot\gamma(t))|\,dt
\\
& = 2 \limsup_{i\to \infty} \int_0^1 |\langle\phi_i,\dot\gamma(t)\rangle|\,dt
\\
& = 2 \int_0^1 |\langle\phi ,\dot\gamma(t)\rangle| \,dt.
\end{align*}
Recall that $\phi$ is locally given by $\sigma\otimes v$, where $v$ is a harmonic $1$--form such that locally $\bv=\pm v$. Since $\dot \gamma$ is tangent to $\ker \bv$, we have $|\langle \phi,\dot\gamma(t)\rangle|=0$ for all $t\in [0,1]$. This implies $d_\infty(x,y)=0$. As a result, every leaf of $\ker \bv$ is contained in a leaf of $\ker \nabla u$ as foliations on $\mathcal{U}$. Therefore, $\ker \bv = \ker \nabla u $ on $\mathcal{U}$. 

Next, we show that $\nabla u$ and $\bv$ differ only by a locally constant multiplication on $\mathcal{U}$. Let $B_r(p)$ be as above. Since $B_r(p)$ is simply connected, both $\nabla u$ and $\bv$ can be lifted to single-valued harmonic $1$--forms.  Locally, write $\bv$ as $\pm v$ and $\nabla u$ as $\pm w$. Both $v$ and $w$ are non-vanishing. Since $\ker \bv = \ker \nabla u$, there exists a non-zero function $g$ such that $v = g\cdot w$. Then we have 
$$
	0 = dv = d(gw) = dg\wedge w + gdw = dg\wedge w,
$$
$$
	0 = d(*v) = d(*gw) = dg\wedge(*w) + g d(*w) = dg\wedge (*w).
$$
As a result, $dg\wedge w = 0$, $dg\wedge (*w)=0$, so $dg=0$, thus $g$ is a constant function. In conclusion, $\nabla u$ and $\bv$ differ only by a locally constant multiplication on $\mathcal{U}$.

Since the zero loci of $\bv$ and $|\nabla u|$ have Hausdorff codimension $2$, the above result implies that $\bv$ is the $\ZT$ harmonic $1$--form associated with $u$ up to a constant multiplication. By Corollary \ref{cor_MH(ell)_distance_function}, we conclude that $\bv = \MH([\ell])$. \qed



\subsection{Proofs of Theorems \ref{thm_map_on_closure} and \ref{thm_factoriztion_tree_maps}}

Define
$$
		\BXi: \BcX\to \BMMc
$$
such that $\BXi$ is given by the (inverse of the) Riemann--Hilbert map on $\cX(\Gamma)$ and  by $\MH$ on $\partial \BcX$. 

\begin{lemma}
	Both $\BcX$ and $\BMMc$ are metrizable.
\end{lemma}

\begin{proof}
	By Corollary \ref{cor_mbM_cpt_Haus}, the space $\BMMc$ is compact and Hausdorff, and hence it is regular. The space $\MMZ$ defined in Section \ref{subsec_ZT_1_form} is separable with respect to the $\MC^0$--topology, because it is a closed subset of the space of continuous $2$--valued sections of $T^*M$, which is a separable metric space. Hence the topology on $\mathbb{M}$ defined in Section \ref{subsec_BMMc} is second countable. By the Urysohn metrization theorem, we conclude that $\BMMc$ is metrizable.

	Similarly, by Lemma \ref{lem_MS_cpt_Haus}, the space $\BcX$  is compact and Hausdorff. Since $C$ is countable, the space $\mathbb{P}^C$ is second countable.  This implies $\BcX$ is second countable. Hence  $\BcX$  is metrizable. 
\end{proof}

The following statement is the first part of Theorem \ref{thm_map_on_closure}:
\begin{theorem}
The map $\BXi$ is continuous. 
\end{theorem}

\begin{proof}
The result follows from Proposition \ref{prop_compare_two_limits} and a formal argument. 
	Since both $\BcX$ and $\BMMc$ are metrizable, 
	we only need to show that if $\{p_i\}$ is a sequence that converges to $p$ in $\BcX$, then $\BXi(p_i)$ converges to $\BXi(p)$. Since $\cX(\Gamma)$ is an open subset of $\BcX$ and the map $\BXi$ is already known to be continuous on $\cX(\Gamma)$, we only need to prove the statement when $p\in \partial \BcX$. 
	We only need to consider two cases: (1) all $p_i$'s are in $\cX(\Gamma)$; (2) all $p_i$'s are in $\partial \BcX$. 
	
	If all $p_i$'s are in $\cX(\Gamma)$, then by the compactness of $\BMMc$, every subsequence of $\BXi(p_i)$ has a convergent subsequence.  By Proposition \ref{prop_compare_two_limits}, this subsequence must converge to $\BXi(p)$. Hence $\BXi(p_i)$ converge to $\BXi(p)$.

	If all $p_i$'s are in $\partial \BcX$, we use an argument by contradiction. Let $d_1$ be a metric for the topology of $\BcX$  and let $d_2$ be a metric for $\BMMc$. Assume there exists $\epsilon>0$ such that there is a subsequence of $\{p_i\}$, which we will still denote by $\{p_i\}$, such that  $d_2(\BXi(p_i), \BXi(p))\ge \epsilon$ for all $i$. For each $i$, since $p_i\in \partial \BcX$, there exists a sequence $\{q_{i,j}\}_{j\ge 1}$ in $\cX(\Gamma)$, such that $\{q_{i,j}\}$ converges to $p_i$ as $j\to\infty$. By Proposition \ref{prop_compare_two_limits}, we have 
	$$
	\lim_{j\to\infty}  \BXi(q_{i,j}) = \BXi(p_i).
	$$
	As a result, we may find a $p_i'$ in the sequence $\{q_{i,j}\}_{j\ge 1}$ such that $p_i'\in \cX(\Gamma)$, $d_1(p_i,p_i')<1/i$, and $d_2(\BXi(p_i'), \BXi(p_i))<\epsilon/2$. The sequence $\{p_i'\}_{i\ge 1}$ then satisfies $\lim_{i\to \infty} p_i' = p$ and $$d_2(\BXi(p_i'), \BXi(p))\ge \epsilon/2$$ for all $i$. This contradicts Proposition \ref{prop_compare_two_limits}. 
\end{proof}

Now we prove the second part of Theorem \ref{thm_map_on_closure}:
\begin{theorem}
	The map $\BXi$ is surjective. 
\end{theorem}

\begin{proof}
	We need to show that every point $q\in \partial \BMMc$ is in the image of $\BXi$. Let $q_i$ be a sequence of points in $\MMc$ such that $q_i\to q$. Let $p_i$ be the preimage of $q_i$ in $\cX(\Gamma)$. By the compactness of $\BcX$, there is a subsequence of $p_i$, which we still denote by $\{p_i\}_{i\ge 1}$, such that $p_i$ converges to a point $p$ in $\BcX$. Since $\BXi$ is continuous, we have $\BXi(p)=q$.
\end{proof}

\begin{proof}
	[Proof of Theorem \ref{thm_factoriztion_tree_maps}]
	By Corollary \ref{cor_MH(ell)_distance_function}, the $\mathbb{Z}/2$ harmonic $1$--form $\bv$ is associated with the harmonic map $u$ up to a non-zero constant multiplication. So Theorem \ref{thm_factoriztion_tree_maps} follows from Theorem \ref{thm_factorize_harm_map_thru_leafs}. 
\end{proof}

\subsection{Non-injectivity of the map $\BXi$}
\label{subsection:non_injectivity}
In general, the map $\BXi$ may not be injective.  
This is certainly well-known, though we have been unable to find an exact
reference.  For completeness we provide the details
of an example. 
In the following, we give a counterexample for the $2$--dimensional analogue.  Multiplying the spaces with $S^1$ will yield a counterexample for the $3$--dimensional case.
Specifically, 
we  construct 
distinct length functions appearing in $\partial \BcX$ with the same image under $\MH = \partial \BXi$.


Let $\Sigma_{g}$ be a Riemann surface of genus $g \geq 3$.
Express the surface as a connect sum: $\Sigma_{g} = \Sigma_{1} \sharp
\Sigma_{g-1}$ of surfaces of genus $1$ and $g-1$, respectively.
This gives an amalgamated product expression for the fundamental group:
$\Gamma=\pi_1(\Sigma_{g}) = \mathbb{F}_2 *_{\mathbb{Z}} \mathbb{F}_{2g-4}$,
where $\mathbb{F}_k$ denotes the free group on $k$ generators. 

Let $\hat\rho_i : \pi_1(\Sigma_{g-1})\to \SL_2(\mathbb{C})$ be a
divergent family of completely reducible representations converging in the
Morgan--Shalen limit to a nonabelian projective length function $\hat\ell$. 
This is possible because $g-1\geq 2$, so we may take
$\hat\rho_i$ to be discrete and faithful, for example. 

Now by the product expression, the $\hat\rho_i$ extend to representations
$\rho_i$ of $\Gamma$, and they clearly converge to a Morgan--Shalen limit
$\ell$, which is the extension of $\hat\ell$. 
Let $\widetilde{\Sigma}_{g}$ be the universal cover of $\Sigma_{g}$, and let $u: \widetilde{\Sigma}_{g} 
\to \MT_{\ell}$ be an equivariant harmonic map with associated $\mathbb{Z}/2$ harmonic
1-form $\bv$. We denote by $\MT_{\bv}$ 
the associated leaf tree defined by $\bv$.
On the other hand, each $\rho_i$ contains the free group $\mathbb{F}_2$ in its kernel, and
therefore the edge stabilizers of $\MT_{\ell}$ all contain free groups. 
By Skora's Theorem \cite{Skora:96},  $\MT_{\ell}$ cannot be dual to a
measured  foliation on $\Sigma_g$. It follows that
$\MT_{\bv}$ and $\MT_{\ell}$ will
necessarily  have different length functions.

Finally, the length function $\ell_{\bv}$ of $\mathcal{T}_\bv$ also appears in the boundary. 
Indeed, as in Example \ref{ex_quadraticdifferential}, $s := \bv \otimes \bv
\in H^0(K^2_{\Sigma_{g}})$ is a holomorphic  quadratic differential on
$\Sigma_g$. The length function $\ell_{\bv}$ is associated to the dual
tree of the vertical measured foliation of $s$, and this appears in the
Thurston boundary of ${\rm (P)SL}_2(\mathbb{R})$ representations of
$\Gamma$. 


Therefore, we have constructed two length functions $[\ell] \neq [\ell_{\bv}]$ such that $\MH([\ell]) = \MH([\ell_{\bv}])$, which implies $\BXi$ is not injective.

\section{Applications }
\label{sec_applications}
In this section, we give several applications of Theorems \ref{thm_map_on_closure} and \ref{thm_factoriztion_tree_maps} and discuss some related results.  

Section \ref{subsec_existenceHakenmanifold} proves Theorem \ref{thm_pi1_inj_Z2_form}, which is an existence result for singular $\mathbb{Z}/2$ harmonic forms on a large class of rational homology spheres.
In Section \ref{subsec_proj_leaf_space}, we show that the projection map
from a simply-connected manifold to the leaf space of a $\ZT$ harmonic
$1$--form may not be harmonic.  In Section
\ref{subsec_ZT_form_simply_connected}, we discuss a folklore conjecture on
the non-existence of $\ZT$ harmonic $1$--forms on $S^3$, and we prove a result about $\ZT$ harmonic forms on simply connected closed manifolds. In Section \ref{subsec_KS_limit_high_reg}, we prove a $\MC^\infty$ convergence result for Korevaar--Schoen limits using results from gauge theory. In Section \ref{subsec_Hub_Mas_map}, we discuss the relationship of Morgan--Shalen limits, singular measured foliations, and the Hubbard--Masur construction. 

\subsection{Existence of singular $\mathbb{Z}/2$ harmonic 1-forms}
\label{subsec_existenceHakenmanifold}

We prove Theorem \ref{thm_pi1_inj_Z2_form} using Culler--Shalen's construction of dual trees associated with essential surfaces.

Let $M$ be a closed oriented 3-manifold with $\Gamma := \pi_1(M)$. Assume $S\subset M$ is an embedded two-sided surface.  
We briefly review the construction of the dual tree of $S$ as in
\cite[Sec.\ 1.4]{shalen02representations}. Let $p: \widetilde{M} \to M$ be
the universal cover and $\widetilde{S} := p^{-1}(S)$ the preimage of $S$.
There is a canonical 1-dimensional simplicial complex $\mathcal{T}_S$
associated with $S$, defined as follows: the vertices are the connected
components of $\widetilde{M} \setminus \widetilde{S}$, and the edges are
the connected components of $\widetilde{S}$. An edge is incidental to a
vertex if the corresponding component of  $\widetilde{S}$ is contained in
the boundary of the corresponding component of $\widetilde{M} \setminus
\widetilde{S}$. The space $\mathcal{T}_S$ is a simplicial tree \cite[Sec.\ 1.4]{shalen02representations}, and the fundamental group $\Gamma$ acts on $\mathcal{T}_S$ by simplicial homeomorphisms. In the following lemma, we show that under the assumptions in Theorem \ref{thm_pi1_inj_Z2_form}, the action of $\Gamma$ on $\mathcal{T}_S$ has no fixed points (see also \cite[Proposition 1.5.2]{shalen02representations}).

\begin{lemma}
	\label{lem_no_fixed_points_pi_inj}
	Suppose $M$ is a rational homology sphere. 
	Suppose $S$ is a closed connected embedded surface in $M$ that is two-sided, $\pi_1$--injective, and does not bound an embedded ball in $M$.  Then the action of $\Gamma$ on $\mathcal{T}_S$ has no global fixed points.
\end{lemma}

\begin{proof}
	Since $M$ is a rational sphere and $S$ is two-sided, $M\setminus S$ must be disconnected. In fact, if $M\setminus S$ is connected, then there is a circle in $M$ that intersects $S$ transversely at one point, which implies that the homology class of $S$ must be non-torsion, contradicting the assumptions. Since $S$ is connected, $M\setminus S$ has two connected components.

	The action of $\Gamma$ on $\mathcal{T}_S$ has a globally fixed vertex if and only if there is a component $C$ of $M\setminus S$ such that $\pi_1(C)\to \pi_1(M)$ is surjective. By the Seifert--van Kampen theorem, this is equivalent to the surjectivity of $\pi_1(S) \to \pi_1(M\setminus C)$.  Since $S$ is $\pi_1$--injective, this condition implies $\pi_1(S)\to \pi_1(M\setminus C)$ is an isomorphism. If $S$ is a sphere, then $M\setminus C$ must be a $3$--ball, which contradicts the assumptions.  If the genus of $S$ is positive, by a standard result in 3-manifold topology (see, for example, \cite[Lemma 3.5]{hatcher2007notes}), the rank of the kernel of $H_1(S)\to H_1(M\setminus C)$ equals the genus of $S$, so $\pi_1(S)\to \pi_1(M\setminus C)$ cannot be isomorphic.
	
	The action of  $\Gamma$ on $\mathcal{T}_S$ has a globally fixed edge if and only if $\pi_1(S)\to \pi_1(M)$ is surjective.  By the Seifert--van Kampen theorem and the assumption that $S$ is $\pi_1$--injective, this implies $\pi_1(S)\to \pi_1(M\setminus C)$ is an isomorphism for each component $C$ of $M\setminus S$.  Hence we get a contradiction by the same argument as before.
\end{proof}


%
\begin{proof}[Proof of Theorem \ref{thm_pi1_inj_Z2_form}]
	Let $\Gamma = \pi_1(M)$. Since $\mathcal{T}_S$ is a simplicial tree, it is complete as a metric space. 
By Lemma \ref{lem_no_fixed_points_pi_inj}, the set of projective length
functions $\MPL(\Gamma)$ is non-empty.  Hence the construction of the map
$\MH: \MPL(\Gamma) \to \MMZ$ in Section
\ref{subsection_construction_boundary_map} implies that $\MMZ$ is
non-empty.  Since $M$ is a rational homology sphere, every non-zero $\mathbb{Z}/2$
harmonic form on $M$ must be non-trivial.
\end{proof}

We remark that essential surfaces also played an important role in the recent construction of new Casson-type invariants by Dunfield and Rasmussen \cite{dunfield2024unified}; see also \cite{dunfield22counting} for results on the counting of essential surfaces.

\subsection{The quotient map to the leaf space}
\label{subsec_proj_leaf_space}
Let $M$ be a closed manifold, $\bv$ be a $\ZT$ harmonic 1-form over $M$,
let $\pi:\tM\to M$ be the universal covering map with $\tbv:=\pi^{*}\bv$,
then by Theorem \ref{thm_pseudo-metric-R-tree}, the leaf space
$(\MT_{\tM,\tbv}, d_{\tM,\tbv})$ is an $\mathbb{R}$--tree.  When $M$ is
2-dimensional, it is proved in \cite[Prop.\
3.1]{{wolf1995harmonic}} that the quotient map from $\widetilde{M}$ to
$\MT_{\tM,\tbv}$ is harmonic. It turns out that this statement does not
hold generally in higher dimensions. The following is a counterexample when $M$ is four dimensional.

\begin{example}
    In \cite[Sec.\ 3, example before Thm.\ 3.2]{bogomolov2011symmetric}, an example of a compact, simply-connected, smooth projective variety $M$ is constructed, such that there is a non-trivial holomorphic section $s$ of $\Omega^1_{M}\otimes \Omega^1_{M}$ that locally has the form $f\mu^2$ with $f$ a holomorphic function and $\mu$ a holomorphic $1$--form. 
    Therefore, the real part of the square-root of $s$ defines a $\mathbb{Z}/2$ harmonic form on $M$, which we denote by $\bv$. The universal cover of $M$ is the same as $M$. Let $\mathcal{T}_{M,\bv}$ be the leaf space of $\bv$ on $M$ as given by Theorem \ref{thm_pseudo-metric-R-tree}.

One can choose the parameters in the construction so that the leaf space of $\ker \bv$ contains infinitely many points. 
To explain this, we need to recall the construction of $M$ from \cite{bogomolov2011symmetric}. Let $A^3$ be a $3$--dimensional abelian variety and let $\hat M$ be a smooth hypersurface in $A^3$. 
Let $\theta = -\id$ be the natural involution on $A^3$, and assume $\hat M$ passes through exactly one fixed point $p$ of $\theta$.  Then $M$ is defined to be the minimal resolution of $\hat M/\theta$. 
Up to multiplication by constants, there is a unique holomorphic $1$--form $w$ on $A^3$ such that $w|_{T_p(\hat M)}=0$. The section $s$ is defined as an extension of the push-forward of $w^2$ to $(\hat M\setminus\{p\})/\theta$.

As a result, the $\mathbb{Z}/2$ harmonic $1$--form $\bv$ is equal to the push-forward of the real part of $\pm w$ on $(\hat M\setminus \{p\})/\theta$. 
Identifying $A^3$ with a quotient of $\mathbb{C}^3$ by a discrete group of translations, one may choose $M$ so that the  space $T_pM$ is given by a linear equation with coefficients in $\mathbb{Q}[\sqrt{-1}]$.  In this case, the kernel of the real part of $\omega$ defines a foliation on $A^3$ such that every leaf is closed. If $L_1,L_2\subset A^3$ are two closed leaves whose images in $A^3/\theta$ are disjoint, are both disjoint from $p$, and both intersect $\hat M$ non-trivially, then the images of $L_i\cap \hat M$ in $\mathcal{T}_{{M},\bv}$ are two distinct points. As a result, $\mathcal{T}_{{M},\bv}$ contains infinitely many distinct points. 

Since ${M}$ is compact, by the maximum principle for harmonic maps, every harmonic map from ${M}$  to an $\mathbb{R}$--tree must be constant. Therefore, the quotient map from ${M}$ to $\mathcal{T}_{{M},\bv}$ cannot be harmonic. \qed
\end{example}

\subsection{$\ZT$ harmonic forms on closed simply connected manifolds}
\label{subsec_ZT_form_simply_connected}
In this subsection, we prove Theorem \ref{thm_ZT_s3}, which states that $\ZT$ harmonic $1$--forms cannot exist on $(S^3,g)$ under certain additional conditions.  In fact, we prove a more general result that holds for closed simply connected manifolds in all dimensions, as stated in Theorem \ref{thm_ZT_sim_conn} below. 

Suppose $M$ is a simply connected closed Riemannian manifold (not necessarily in dimension $3$), and $\bv$ is a $\mathbb{Z}/2$ harmonic $1$--form on $M$ with zero locus $Z$. Let $\mu_\bv$ be the transverse measure on $M\setminus Z$ defined by $\bv$ (see \eqref{eqn_transverse_measure_from_bv}), and let $d_{M,\bv}$ be the corresponding pseudo-metric on $M$ (see \eqref{eq_lengthfunction_defined_by_Z2harmonicform}); we also use $d_{M,\bv}$ to denote the induced metric on the leaf space $\MT_{M,\bv}$. 

\begin{definition}
	\label{def_cylindrical_zero}
	We say that a zero point $p\in Z$ of $\bv$ is \emph{cylindrical}, if there exists a local coordinate chart $(x_1,\dots,x_n)$ on a neighborhood of $p$ such that $\bv=\varphi \cdot Re(d\sqrt{z^{k}})$ for some integer $k\ge 3$, where $\varphi$ is a smooth and non-vanishing function, and $z=x_1+x_2\sqrt{-1} \in\mathbb{C}$.
\end{definition}

The following theorem is a generalization of Theorem \ref{thm_ZT_s3}.
\begin{theorem}\label{thm_ZT_sim_conn}
	Suppose $M$ is a closed, simply-connected, Riemannian manifold. Then there is no $\ZT$ harmonic $1$--form $\bv$ on $M$ such that the following conditions hold at the same time.
	\begin{enumerate}
		\item Every point $p$ of $\bv$ is cylindrical.
		\item For every arc $\gamma$ in $M\setminus Z$ transverse to $\ker \bv$, where $Z$ denotes the zero locus of $\bv$, we have 
		$$
		\mu_\bv(\gamma) = d_{M,\bv}(x,y).
		$$
	\end{enumerate}
\end{theorem}

\begin{remark}
	Condition (1) implies that the hypotheses of Theorem \ref{thm_pseudo-metric-R-tree} hold.
	By Lemma \ref{lem_compare_measure_and_image}, Condition (2) always holds if the $\ZT$ harmonic form is associated with an equivariant harmonic map from a universal cover of a closed $3$--manifold to a tree.
\end{remark}

Theorem \ref{thm_ZT_sim_conn} will be a straightforward consequence of the following lemma.

\begin{lemma}
	\label{lem_sufficient_cond_harm}
	Assume $M$ is a simply connected (but not necessarily closed) manifold and $\bv$ satisfies the conditions in Theorem \ref{thm_ZT_sim_conn}. Let $\pi_\bv$ be the projection map from $M$ to the leaf tree $\MT_{M,\bv}$. Then germs of convex functions on $\MT_{M,\bv}$ pull back to germs of subharmonic functions on $M$.
\end{lemma}

Here, a function $f$ on an $\mathbb{R}$--tree $\MT$ is called \emph{convex}, if for every geodesic $\gamma:[a,b]\to \MT$ parametrized by arc length and every $\lambda\in[0,1]$, we have $$\lambda f\circ\gamma(a) + (1-\lambda) f\circ\gamma(b) \ge  f\circ\gamma(\lambda a+(1-\lambda)b).$$

\begin{remark}
The condition that germs of convex functions pull back to germs of
subharmonic functions is a local condition that can be verified near
every point. It is also well-known that for maps between
manifolds, this condition is equivalent to the map being harmonic (see,
for example, \cite[Thm.\ 3.4]{ishihara1979mapping}). 
\end{remark}

\begin{proof}
	By Condition (2) in Theorem \ref{thm_ZT_sim_conn}, the map $\pi_\bv$ is regular and harmonic near every non-zero point of $\bv$, so it satisfies the desired properties on the complement of the zero locus of $\bv$. 
	
	Let $p$ be a zero point of $\bv$. Take a coordinate chart $(x_1,x_2,x_3,\dots,x_n)$ containing $p$ such that Definition \ref{def_cylindrical_zero} holds. Let $k$ be the integer in Condition (1). Let $U_\epsilon$ be the open neighborhood of $p$ given by $|x_i|<\epsilon$ $(i=1,\dots,n)$ in this chart. 
	
	Let $\MT_{U_\epsilon,\bv}$ denote the leaf space of $\bv$ on $U_\epsilon$, then  $\MT_{U_\epsilon,\bv}$ is given by the union of $k$ segments with one end point identified.  The identified end point is the image of $p$, which we denote by $\hat p$. 
	In the following, we will call each segment in $\MT_{U_\epsilon,\bv}$ a \emph{branch}. The $k$ branches admit a natural cyclic order.
	
	Let $\MT_\epsilon$ denote the image of $\MT_{U_\epsilon,\bv}$ in $\MT_{M,\bv}$. We will abuse notation and use $\hat p$ to also denote the image of $p$ in $\MT_{M,\bv}$. 
	
	By Condition (2), for $\epsilon$ sufficiently small, the map of each branch of $\MT_{U_\epsilon,\bv}$ into $\MT_{\epsilon}$ is an isometric embedding.  Therefore, the map from $\MT_{U_\epsilon,\bv}$ to $\MT_{\epsilon}$ is a quotient map that identifies pairs or groups of branches.
	Moreover, Condition (2) also implies that neighboring branches of $\MT_{U_\epsilon,\bv}$ cannot be identified in $\MT_{\epsilon}$.

	Let $f$ be a convex function defined on a neighborhood of the closure of $\MT_\epsilon$, we show that $f$ pulls back to a subharmonic function near $p$.  
	
	We first consider the case that $f(x)=- \mathrm{dist}_{\MT_\epsilon}(x,\hat{p})$ on one of the branches of $\MT_{\epsilon}$ (which we denote by $I$), and $f(x)=\mathrm{dist}_{\MT_\epsilon}(x,\hat{p})$ on $\MT_\epsilon\setminus I$.
	 
	  Let $\tilde f$ be the pull-back of $f$ to $U_\epsilon$. Then the pre-image of $\hat p\in \MT_{U_\epsilon,\bv}$ in $M$ cuts $U_\epsilon$ into $k$ domains, which we denote by $\Omega_1,\dots,\Omega_k$ in cyclic order. 
	   Then there exists a set $S\subset \{1,\dots,k\}$, such that 
	   \begin{equation}
	   	\label{eqn_tilde_f}
	   		   \tilde{f}(x) =
	   		   \begin{cases}
	   		     -d_{U_\epsilon,\bv}(x,p) \quad \text{ if } x\in \Omega_i, \,\, i\in S \\
	   		       d_{U_\epsilon,\bv}(x,p) \quad \text{ if } x\in \Omega_i, \,\, i\notin S.
	   		    \end{cases}
	   \end{equation}
	   In the following, the subscripts will always be interpreted modulo $k$.
Since neighboring branches of $\MT_{U_\epsilon,\bv}$ have distinct images in $\MT_{\epsilon}$, we have $\{i,i+1\}\not\subseteq S$ for every $i$.
	
	By \cite{ishii1995equivalence}, we only need to show that 
	\begin{equation}\label{eqn_tilde_f_boundary_int_non_negative}
		\int_{\partial U_\epsilon} \langle \nabla \tilde{f}, n \rangle \ge 0,
	\end{equation}
	for all $\epsilon$,
	where $n$ denotes the outward unit normal vector field on $\partial U_\epsilon$. For each $\Omega_i$, we decompose $\partial \Omega_i$ into $\partial \Omega_i = F_i \cup F_i' \cup F_i''$, where $F_i = \partial \Omega_i \cap \partial U_\epsilon$, $F_i'= \Omega_i\cap\Omega_{i-1}$, $F_i''=\Omega_i\cap\Omega_{i+1}$.  It is clear that $F_i'=F_{i-1}''$ for all $i$. 
	
	We have 
	$$
	\int_{\partial \Omega} \langle \nabla\tilde{f}, n\rangle  = \sum_{i=1}^n 	\int_{F_i} \langle \nabla\tilde{f}, n\rangle .
	$$
	Define 
	$$
	E_i = \int_{F_i'} |\langle\bv,n\rangle|,
	$$
	where $n$ is a unit normal vector field of $F_i'$.
	By \eqref{eqn_tilde_f} and the fact that $d_{U_\epsilon,\bv}(x,p)$ is a harmonic function on $x$ in the interior of each $\Omega_i$, we have 
	$$
	\int_{F_i} \langle \nabla \tilde{f}, n\rangle= 
	\begin{cases}
		-E_i - E_{i+1} \quad \text{ if }  i\in S, \\
		E_i + E_{i+1} \quad \text{ if }  i\notin S. \\
	\end{cases}
	$$
	Hence \eqref{eqn_tilde_f_boundary_int_non_negative} follows from the fact that $\{i,i+1\}\not\subseteq S$ for all $i$.
	
	This proves that the pull-back of $f$ is subharmonic near $p$ when $f$ has the special form given as above. In general, assume $g$ is an arbitrary convex function that is defined on a neighborhood of the closure of $\MT_\epsilon$.  Then there exist constants $c_1\ge 0$, $c_2$, and a function $f$ having the form above, such that $g\ge c_1f+c_2$ on $\MT_\epsilon$ and $g(\hat p) = c_1f(\hat p)+c_2$. Let $\tilde{g}$ denote the pull back of $g$ to $U_\epsilon$. Let $h_g$ denote the harmonic function on $U_\epsilon$ such that $h_g=\tilde g$ on $\partial U_\epsilon$, let $h_f$ denote the harmonic function on $U_\epsilon$ such that $h_f=\tilde f$ on $\partial U_\epsilon$. Then we have
	$$
	\tilde g(p) = c_1 \tilde f(p) + c_2 \le c_1 h_f(p)+ c_2 \le h_g(p).
	$$
	Since the above inequality holds for all $p\in |\bv|^{-1}(0)$ and $\epsilon$, we conclude that $\tilde g$ is subharmonic, and the lemma is proved.
\end{proof}

\begin{proof}[Proof of Theorem \ref{thm_ZT_sim_conn}]
	Assume $\bv$ is a $\ZT$ harmonic form satisfying the conditions of Theorem \ref{thm_ZT_sim_conn}. Then for every $q\in \MT_{M,\bv}$, the distance to $q$ is a convex function on $\MT_{M,\bv}$. By Lemma \ref{lem_sufficient_cond_harm} and the maximum principal, the distance function to $q$ must be constant on the image of $\pi_\bv$.  Since this holds for all $q$, the map $\pi_\bv$ must be constant and hence $\MT_{M,\bv}$ contains only one point. By Condition (2) in Theorem \ref{thm_ZT_sim_conn}, the leaf space $\MT_{M,\bv}$ contains infinitely many points, which yields a contradiction. 
\end{proof}

Theorem \ref{thm_ZT_s3} then follows as a special case of Theorem \ref{thm_ZT_sim_conn}.

\subsection{$\MC^\infty$--convergence for Korevaar--Schoen limits}
\label{subsec_KS_limit_high_reg}
Theorem \ref{thm_map_on_closure} establishes a bridge between the analytic and algebraic compactifications of the moduli space of flat $\SLC$ connections on $3$--manifolds. This connection allows us to prove new results on one side using results from the other.
Let $M$ be a closed $3$--manifold, let $\Gamma=\pi_1(M)$, and let $\{p_i\}$ be a sequence of points in $\cX(\Gamma)$ that converges to a point in $\partial \overline{\cX(\Gamma)}$.
Let $u_i: \widetilde{M} \to \mathbb{H}^3$ be the corresponding equivariant
harmonic maps (see the discussion in Section \ref{subsubsec_MS_lim}). Let
$E(u_i)$ be the energy of $u_i$ on a fundamental domain. By  \cite[Thm.\
2.2]{daskalopoulos1998character} (see also Section \ref{subsubsec_MS_lim}),
after rescaling, the limit of $u_i$ defines an equivariant harmonic map $u:
\widetilde{M} \to \mathcal{T}$ to a tree.  By \cite[Thm.\
3.9]{korevaari1997global}, for every smooth vector field $W$ on
$\widetilde{M}$, the sequence $|\nabla u_i(W)|/E(u_i)^{1/2}$ converges
weakly to $|\nabla u(W)|$ in $L^2$. We use a result of Parker \cite[Thm.\  1.3]{parker2023concentrating} to show that, in fact, $|\nabla u_i(W)|/E(u_i)^{1/2}$ also converges to $|\nabla u(W)|$ in $C^\infty$ on compact subsets of the complement of the zero locus of $|\nabla u|$.

\begin{corollary}
	\label{cor_C^infty_KS_lim}
	Let $u_i, u$ be as above. Let $\mathcal{U}$ be the open subset of $\widetilde{M}$ where $|\nabla u| \neq 0$. Then for every smooth vector field $W$, we have $|\nabla u_i(W)|/E(u_i)^{1/2}$ converges in $\MC^\infty$ to $|\nabla u(W)|$ on compact subsets of $\mathcal{U}$. 
\end{corollary}

\begin{proof}
	Without loss of generality, assume $K$ is a compact set contained in a fundamental domain of $\widetilde{M}$, and assume that $W$ is $\pi_1(M)$ invariant.  Then $|\nabla u_i(W)|$ and $|\nabla u(W)|$ reduce to functions on $M$, and $W$ reduces to a vector field on $M$. Assume $(A_i,\phi_i)$ is the solution to \eqref{eqn_KW_system} corresponding to $u_i$, and let $\bv\in \MMc$ be the limit of the sequence $(A_i, \phi_i)$ in $\partial \BMMc$. By Theorem \ref{thm_map_on_closure} and Corollary \ref{cor_MH(ell)_distance_function}, we have $|\bv(W)| = |\nabla u(W)| /\|\nabla u\|_{L^2(M)}$. 
	
	We abuse notation and also use $K$ to denote the image of $K$ in $M$. 
	By \cite[Thm.\ 1.3]{parker2023concentrating}, we know that every subsequence of $|\nabla \phi_i(W)|/\|\phi_i\|_{L^2(M)}$ has a subsequence that converges in $\MC^\infty$ on $K$. Hence, by \eqref{eqn_relation_u_phi}, every subsequence of $|\nabla u_i(W)|/E(u_i)^{1/2}$ has a subsequence that converges in $\MC^\infty$ on $K$. Since $|\nabla u_i(W)|/E(u_i)^{1/2}$ weakly converges to $|\nabla u(W)|$, the $\MC^\infty$ limit of the subsequence must also be $|\nabla u(W)|$. Therefore, $|\nabla u_i(W)|/E(u_i)^{1/2}$ converges to $|\nabla u(W)|$ in $\MC^\infty$ on $K$.
\end{proof}

\subsection{Properties of the  boundary of the analytic compactification $\partial \BMMc$}
\label{subsec_Hub_Mas_map}
In this subsection, we discuss several properties of the  $\ZT$ harmonic
1-forms that actually appear  in $\partial \BMMc$. We 
compare them with classical results of measured foliations and the Hubbard-Masur map \cite{hubbard1979quadratic} over Riemann surfaces. 
\subsubsection{Leaf spaces of the analytic boundary}
Observe that by Theorem \ref{thm_factoriztion_tree_maps}, if a
$\mathbb{Z}/2$ harmonic form $\bv$ arises as a limit in Taubes'
compactification, then  the leaf space of $\bv$ has similar properties to the leaf spaces of quadratic differentials.  
\begin{corollary}
	\label{cor_condition_MMZ_taubes_lim}
	Suppose $\bv\in \partial \BMMc$ and let $Z$ be the zero locus of $\bv$.
    Let $\widetilde{M}$ be the universal cover of $M$ and let $\tilde{\bv}$
    denote the pull-back of $\bv$ to $\widetilde{M}$. Let $\mu_{\tilde \bv}$, $d_{\widetilde{M},\tilde \bv}$ be defined as in \eqref{eqn_transverse_measure_from_bv},  \eqref{eq_lengthfunction_defined_by_Z2harmonicform},
    and let $\MT_{\widetilde{M}, \tilde\bv}$  be the leaf space given by Theorem \ref{thm_pseudo-metric-R-tree}. Then 
	\begin{enumerate}
		\item the quotient map from $\widetilde{M}$ to $\MT_{\widetilde{M},
            \tilde\bv}$ is harmonic; and
		\item  if $\gamma$ is an arc in $\widetilde{M}$ from $x$ to $y$  that is transverse to $\ker \tilde\bv$, then $\mu_{\tilde \bv}(\gamma) = d_{\widetilde{M}, \tilde\bv}(x,y)$.
	\end{enumerate}
\end{corollary}

\begin{proof}
Statement (1) follows from Theorem \ref{thm_factoriztion_tree_maps}. By Theorem \ref{thm_map_on_closure}, $\bv$ is in the image of $\MH$. Statement (2) then follows from Corollary \ref{cor_MH(ell)_distance_function} and Lemma \ref{lem_compare_measure_and_image}.
\end{proof}

\subsubsection{Canonical measured foliations and the Hubbard--Masur map}
We give a geometric interpretation for the restriction of $\BXi: \BcX \to \BMMc$ to the boundary $\partial \BcX$, which equals $\MH|_{\partial\BcX}$.
In Morgan--Shalen's theory \cite{morganshalen1984valuationstrees,
morganshalen88degeneration, morganshalen88degernerationtwo}, if $[\ell] \in
\partial\BcX$ is in the limit with the
corresponding $\Gamma$--tree $\MT_{\ell}$, the transverse equivalence map construction
in \cite[Sec.\ I.2]{morganshalen88degernerationtwo} defines a measured lamination from the harmonic map.

By comparison, for each $[\ell] \in \partial\BcX$, the corresponding $\ZT$ harmonic form $\MH([\ell])$ defines a (singular) measured foliation on $\widetilde{M}$ as in Section \ref{subsec_ZT_form_to_foliation}. Therefore, the map $\MH$  constructed in Section \ref{subsection_construction_boundary_map} canonically associates a (singular) measured \emph{foliation} with every Morgan--Shalen limit point.

In fact, the boundary map $\MH$ can be viewed as a generalization of the Hubbard--Masur map. The original Hubbard--Masur map \cite{hubbard1979quadratic} establishes a homeomorphism between the space of projective classes of measured foliations (up to equivalence) and the space of quadratic differentials on a Riemann surface $\Sigma$. As shown in Corollary \ref{cor_MH(ell)_distance_function}, the map $\MH$ generalizes this concept to three dimensions, where harmonic maps to $\mathbb{R}$--trees take the place of measured foliations.

\bibliography{references}
\bibliographystyle{amsalpha}
\end{document}